\documentclass[11pt]{article}
\usepackage{amsmath, amscd,amsbsy, amsthm, amsfonts, amssymb}
\usepackage{rotating}
\usepackage{dcolumn,multicol,url}

 \theoremstyle{definition}

 \numberwithin{equation}{section}

 \newcommand{\A}{\mathcal{A}}

\newcommand{\Gsp}{G_n^{spin}(X)}
\newcommand{\Gp}{G_n^{pin}(X)}

\def\C{{\mathbf C}}

\def\Q{{\mathbb Q}}
\def\R{{\mathbb R}}

\def\Z{{\mathbb Z}}

\def\Q{\mathbb{Q}}

\def\be{\begin{equation}}
\def\ee{\end{equation}}
\def\ba{\begin{eqneqnarray}}
\def\ea{\end{eqneqnarray}}
\def\tilde{\widetilde}

\def\e1{\epsilon}
\def\AAl{\mathcal{A}_{\lambda}}
\def\A0{\stackrel{\circ}{\AAl}}

\def\o1{\widetilde{\Omega}}
\def\01{\widetilde{\Omega}}
\def\c1{\gamma}

\def\g1{\Sigma}

\def\l1{\Lambda}

\def\v1{\varphi}

\def\d1{\delta}

\def\f1{\frac}
\def\t1{\theta}
\def\b1{\beta}
\def\bar{\overline}
\def\bs{\begin{eqneqnarray*}}
\def\es{\end{eqneqnarray*}}
\def\m1{\Theta}
\def\w1{\wedge}

\begin{document}

\title{Quadratic Functions of Cocycles and Pin Structures}

\author{Greg Brumfiel and John Morgan}

\maketitle

\tableofcontents

\section{Introduction}

\subsection{Preview}

The central result of this paper is a combinatorial characterization of equivalence classes of Pin$^-$ structures on a compact triangulated\footnote{Throughout the paper, all manifolds are compact. By a {\em triangulation}, we mean an ordered simplicial structure.  See \S1.3 for details.  A Pin$^{-}$ structure for us means the Wu class $v_2(M)$ is killed.}  $n$-manifold $M$, possibly with boundary, as certain $\Z/4$-valued `quadratic functions'  $Q$ defined on relative $(n-1)$-cocycles  $Z^{n-1}(M,  \partial M; \Z/2)$.   The function $Q$ should satisfy two conditions involving the $\cup_i$ products of Steenrod.  Specifically, we will prove the following.\\

\noindent CLAIM 1.1(P):  Equivalence classes of Pin$^-$ structures on $M^n$ are in canonical bijective correspondence with functions $$Q\colon Z^{n-1}(M, \partial M; \Z/2) \to  \Z / 4$$ that satisfy $$Q(p+q) = Q(p) + Q(q) + 2\int_{[M, \partial M]}\  p \cup_{n-2} q \in \Z/4$$ and $$Q(dc) = 2\int_{[M, \partial M]}\ Sq^2c\ =  2\int_{[M, \partial M]}\  c \cup_{n-4} c + c \cup_{n-3} dc \in \Z/4.$$

We will refer to such $Q$ as  $\it{quadratic\ functions}$ on $M$.  The integrals take values in $\Z/2$, and the 2 in front of the integrals refers to the inclusion $\Z/2 \to \Z/4$.    In \S2.1 we will verify that the formula for $Q(dc)$ depends only on $dc$, not on $c$. We will also see that the two conditions on $Q$ are consistent for sums $dc + dc'$. \\

When does such a quadratic $Q$ exist?  Taking $c = 0$ in the second condition, we see  $Q(0) = 0$.  Then taking $dc = 0$, that is, $c \in Z^{n-2}(M, \partial M; \Z /2)$, we see $$0 = Q(0) = 2\int_{[M, \partial M]}\ c \cup_{n-4} c = 2\langle Sq^2c, [M, \partial M]\rangle.$$
Thus $M$ must satisfy $v_2(M) = 0$, which is the condition that $M$ admits Pin$^-$ structures and which we now assume.\\

The characteristic class $v_2(M)$ coincides with $w_2$ of the stable normal bundle of $M$, which also coincides with $w_2 + w_1^2$ of the tangent bundle.  We have no reason to get involved with the discussion of Pin$^+$ structures, which exist when $w_2$ of the tangent bundle vanishes.  Therefore, from here on we will simply refer to a Pin$^-$ structure as a Pin structure. \\

If $M^n$ is oriented, we will see that the values of such $Q$ must  lie in $2( \Z/2) \subset \Z/4$. An oriented Pin manifold is a Spin manifold.  In that case we can restate the Claim in a slightly neater form.\\

\noindent CLAIM 1.1(S):  Equivalence classes of Spin structures on an oriented manifold $M^n$ are in canonical bijective correspondence with functions $$Q\colon Z^{n-1}(M, \partial M; \Z/2) \to  \Z / 2$$ that satisfy $$Q(p+q) = Q(p) + Q(q) + \int_{[M, \partial M]}\  p \cup_{n-2} q \in \Z/2$$ and $$Q(dc) = \int_{[M, \partial M]}\ Sq^2c\ =  \int_{[M, \partial M]}\  c \cup_{n-4} c + c \cup_{n-3} dc \in \Z/2.$$

Actually, the terminology of Spin and Pin structures is not so appropriate, since smoothness has nothing to do with these quadratic functions.  Manifold can mean PL manifold, or a more general space like a triangulated homology manifold or a Poincare duality space.  Any structure where killing $v_2$ makes sense is suitable, but we will continue to use the Spin and Pin terminology.\\

If $v_2(M) = 0$ then one can construct quadratic functions on $M$ as follows.  We assume for now the details to be checked in \S2.1 concerning the consistency of the definition of quadratic functions.  Use the second condition to define $Q(dc)$.    If $\{p_j\}$ are cocycles representing a basis of $H^{n-1}(M, \partial M ; \Z /2)$, then the values $Q(p_j)$ will determine all $Q(p)$, using the quadratic condition.  Taking $p = q$ and using the fact that $p \cup_{n-2} p = Sq^1p$, we see that $M$ orientable  actually forces $0 = Q(2p) = 2Q(p)$. The $Q(p_j)$ can then be assigned arbitrarily in $2 (\Z/2) \subset \Z/4$.  In the general Pin case, one can still assign the values $Q(p_j) \in \Z/4$ arbitrarily, subject to the constraint  $2Q(p_j) = 2 \int_{[M, \partial M]}\ Sq^1p_j \in \Z/4$. We will omit the easy but somewhat tedious verification that the construction here does define quadratic functions on $Z^{n-1}(M, \partial M; \Z/2)$.\\

A second  $Q'$ differs from $Q$ by a linear function $H^{n-1}(M, \partial M; \Z/2) \to \Z/2$.  This means $Q'(p) = Q(p) + 2\langle pa, [M, \partial M]\rangle$, for some $a \in H^1(M; \Z/2)$. The arguments here show that the existence and enumeration of quadratic functions of cocycles on $M$ does not depend on the choice of  simplicial structure.  Note that if $Q$ is a quadratic function then so is $-Q$.  In fact $$-Q(p) = Q(p) + 2Q(p) = Q(p) + 2 \int_{[M, \partial M]} Sq^1p = Q(p) + 2\langle pv_1, [M, \partial M]\rangle ,$$ where $v_1 = w_1$ is the first Wu class and Stiefel-Whitney class of $M$.\\

If $\sigma^*$ denotes the cochains that are vector space duals of  single $(n-2)$-simplices $\sigma$ in $M - \partial M$, then given the first quadratic condition on $Q$, the second (somewhat mysterious) condition on $Q(dc)$ is equivalent\footnote{This is an observation of some importance.} to $Q(d\sigma^*) = 0$, all $\sigma$. The proof that $Q(d\sigma^*) = 0$ requires looking at the definition of the $n$-cochains $\sigma^* \cup_{n-4}\sigma^*$ and $\sigma^* \cup_{n-3}d\sigma^*$ as a sum of products of evaluations of $\sigma^*$ and $d\sigma^*$ on faces of an $n$-simplex.  In fact, both these $\cup_i$ products vanish  on all $n$-simplices, so the integrals certainly vanish.\\

For closed surfaces, that is,  $n = 2$ and $\partial M = \emptyset$, the Claims 1.1(P) and 1.1(S) are well-known, [4]. In that dimension the second condition just reduces to $Q(dc) = 0$ and the first condition then says $Q$ is a $\Z/4$-valued quadratic refinement of the (non-singular) cup product pairing on $H^1(M; \Z/2)$.   If $n > 2$ then $Q$ is only defined on cocycles, not on cohomology classes, since possibly $Q(p+ dc) \not= Q(p)$. In all cases except closed surfaces there is no reason the pairing $\int_{[M, \partial M]}\ p \cup_{n-2} q$ should have any non-degenerate properties.\\

We will eventually prove that if $M$ is given a Pin structure then a $\it{canonical}$ such $Q$ exists. Equivalence classes of Pin structures are always a torsor of $H^1(M; \Z/2)$. Thus, the canonical bijective correspondence part of  Claim 1.1(P) amounts to showing that changing the Pin structure on $M$ by a cohomology class changes the canonical quadratic function by the same class.   Our proof clarifies this point in some detail. Because $Q$ and $-Q$ are different in the general Pin case, the canonical quadratic function is not quite as well-defined as in the Spin case.  It is necessary to first make a normalizing choice with $\R P^2$, that is, declare which quadratic function will be assigned to which Pin structure.  But after that choice, there is a canonical quadratic function on any Pin manifold, with or without boundary.\\

It turns out that a proof of  Claim 1.1(S) for Spin manifolds can be given that is quite a bit more direct than the only proof we could come up with for general Pin manifolds. This direct proof in the Spin case is essentially due to Kapustin, [3].  We will  explain briefly in the next few paragraphs the difference between the Spin and Pin cases.\\

In our two papers [1] and [2] studying the Pontrjagin dual of the (reduced) Spin bordism functor $\widetilde{\Omega}^{spin}_n(X)$ in dimensions $n = 3$ and $n = 4$, an important ingredient was describing a certain functorial subgroup of these Pontrjagin duals.  This subgroup is actually well-understood in all dimensions. For a simplicial complex $X$, a group $G_n(X)$ is constructed from certain cocyles and cochains, along with an isomorphism\footnote{\label{circle} We identify $\R/ \Z$  with the  unit circle $S^1$ via complex exponentiation $e^{2\pi it}$.  For our purposes,  $\R/ \Z$ is advantageous for the values of characters as it allows for simultaneous additive notation for cochain groups with values in $\Z, \Z/n$, and $\R / \Z$.}
\begin{equation}
G_n(X) \simeq {\rm Hom}(\widetilde{\Omega}^{spin}_n(X) / Image(\widetilde{\Omega}^{spin}_n(X^{n-2})),\ \R / \Z),
\end{equation}
where $X^{n-2}$ is the $n-2$ skeleton of $X$.  In fact, our work on the Pontrjagin dual of Spin bordism began after the physicist Anton Kapustin asked us if we could extend his understanding of this subgroup   to the full Pontrjagin dual in dimensions 3 and 4. The physics paper [3]  contains, among many other things,  discussions essentially amounting to a combinatorial description of $G_n(X)$ and the map (1.1) to the Pontrjagin dual of Spin bordism.  Some of the steps in the  construction of this map amount to a definition of a canonical quadratic function for closed Spin $n$-manifolds.  The proof that the map is an isomorphism to the indicated subgroup of the Pontrjagin dual of Spin bordism required a second definition of the Kapustin quadratic function on Spin manifolds.\\

We included the details of the combinatorial description of $G_n(X)$, along with the connections with  quadratic functions and Spin structures on closed manifolds, in the Appendix \S10.3 of our $4D$ paper [2].  We reproduce those details in this paper, but we do more.\\

First, it is quite easy to extend that basic discussion to Spin manifolds with boundary $(M, \partial M)$.  We also extend the discussion to Pin manifolds, with or without boundary.  In particular, for a Pin manifold we define a function $Q\colon Z^{n-1}(M, \partial M; \Z/2) \to \Z/ 4$, using an analogue of the second construction for Spin manifolds.  But we cannot directly prove our function $Q$ is quadratic, because the first Kapustin construction  does not seem to have an analogue for Pin manifolds.  Therefore, we need to use some tricky stable homotopy theory to incorporate our function $Q$ as part of a Pin bordism analogue of the isomorphism (1.1) for all spaces $X$.  In fact, we extend (1.1)  to pairs $(X, Y)$. Specializing to a Pin manifold $X = M^n$, the identity relative bordism element $Id\colon (M^n, \partial M^n) \to (M^n, \partial M^n)$  defines by duality a linear function from the Pin version of $G_n(M^n, \partial M^n)$ to $\R / \Z$.  But a linear function on $G_n(M^n, \partial M^n)$ turns out to be essentially the same thing as a quadratic function on $M^n$, so we conclude our function $Q$ is necessarily quadratic.\\ 

The stable homotopy theory needed to carry out this strategy is the identification of classifying spaces for the functors on the right side of the map (1.1), in both the Spin and Pin cases.  We construct these classifying spaces as simplicial sets that are two stage Postnikov towers $E_n$, which together form a spectrum as $n$ changes.  To relate these spaces for different $n$, we exploit a cochain suspension operation $s$ introduced in \S4, that we also used in our paper [2] for somewhat similar purposes.  The crucial property of $s$ is the formula (4.1) of \S4, expressing a  commutativity of $s$ and $\cup_i$ products.\\ 

We also use the cochain suspension operation $s$ to directly explain how quadratic functions on $M$ explicitly determine quadratic functions on $\partial M$, and also on codimension 0 submanifolds of $M$.  Obviously, a classical Pin structure on $M$ determines a Pin structure on $\partial M$, so it is good to understand this in the language of the quadratic functions.  We  give several other examples of direct manipulations of quadratic functions that correspond to equivalent manipulations with Pin structures.\\

One can find in the literature various discussions of combinatorial equivalents to Spin structures and maybe Pin structures, especially in low dimensions.  We find our discussion of Pin structures rather clean.  Given a simplicial structure, a quadratic function is one thing, not an equivalence class of things.  However, there is a price.  The description is in terms of the somewhat mysterious $\cup_i$ cochain  products of Steenrod.  In theory, one could translate everything we say into a more geometric discussion of manipulations with numbers assigned to 1-dimensional graphs in $n$-manifolds, dual to $n-1$  cocycles and coboundaries.  For example, a 1-cycle dual to $d\sigma^*$, where $\sigma$ is an $n-2$ simplex, is a little circle link of $\sigma$ passing through the $n$-simplices that have $\sigma$ as a face. This circle trivially bounds a framed disk, which explains why $Q(d\sigma^*) = 0$ is reasonable.  Locally Pin structures are pretty simple, it is  fitting them together globally that is tricky. Part of the physics paper [3] looks at such geometric constructions in the Spin case to some extent.  But the translation is not easy.  It is better to accept  the multivariable $\cup_i$ cochain operations as an important ingredient of nature to be exploited.  After all, the revered Steenrod squares themselves are defined as $\cup_i$ products at the cocycle level,  but generally are treated as black boxes with certain amazing properties.  If one looks into the actual cochain formulas for Steenrod squares, they are pretty complicated.\\ 

\subsection{Organization of the Sections of the Paper}

We conclude this introductory Chapter in the next section with a few remarks about ordered simplicial complexes and cochain operations.  In the first few sections of Chapter 2 we give a direct proof of the main Claim 1.1(S) for Spin manifolds, by explicitly constructing canonical quadratic functions.  In \S2.5 we define a canonical  function for a Pin manifold, but we cannot yet prove much about it, especially that it is quadratic.\\ 

In Chapter 3 we construct functorial groups $G^{spin}_n(X, Y)$  and $G^{pin}_n(X, Y)$ for all simplicial  pairs $(X, Y)$ from cochains and cocycles.  When $(X, Y) = (M^n, \partial M^n)$ is an $n$-manifold, it is trivial that certain elements in the linear Pontrjagin dual of the  group $G^{spin}_n(M, \partial M)$ correspond to $\Z/2$-valued quadratic functions on $M$, and elements in the  dual of $G^{pin}_n(M, \partial M)$ correspond to $\Z/4$-valued quadratic functions on $M$. This is stated formally in \S3.3.  Then  in \S3.4 we show using singular complexes that the functors $G_n(X, Y)$ are homotopy functors.\\

In Chapter 4 we introduce a natural cochain suspension operation $$s \colon C^k(X) \to C^{k+1}(\Sigma X).$$  This is a classical construction, but we clarify and emphasize how $s$ commutes with the cochain $\cup_i$ operations, a property that has perhaps been somewhat overlooked.\\

In Chapter 5 we review some elementary simplicial theory of Postnikov towers.  This includes a brief discussion of 2-stage Postnikov tower simplicial $H$-spaces.  Then we explain how the cochain suspension operation $s$ of Chapter 4 has a more homotopy theoretic interpretation that relates a Postnikov tower to its loop space.  One motivation for this is that the homotopy theory adds some clarity to our later discussion in Chapter 7 that uses cochain suspension to relate quadratic functions on manifolds to quadratic functions on their boundaries and on codimension zero submanifolds.\\

In Chapter 6 we define and study the 2-stage Postnikov towers $E^{spin}_n$ and $E^{pin}_n$ that classify the homotopy functors $G^{spin}_n(X)$  and $G^{pin}_n(X)$ of Chapter 3.  These towers are simplicial $H$-spaces. But the tricky thing is that we know {\em in advance} that these spaces also represent the Pontrjagin duals of the quotients of Spin and Pin bordism occurring on the right in the map (1.1).  Therefore, the groups $G^{spin}_n(X)$  and $G^{pin}_n(X)$ are also naturally isomorphic to these Pontrjagin duals, by an essentially unique isomorphism.  In the Spin case, we can directly define the map (1.1) and prove it is a group isomorphism.  In the Pin case, we only know there is an isomorphism, but we know just enough about it to conclude that the isomorphism must be essentially given by the would-be quadratic function $Q$ defined in \S2.5, and therefore that function is necessarily quadratic. The  details of the proof of the main result, Claim 1.1(P), concerning quadratic functions and Pin structures are given in \S\S6.4-6.7.\\

Finally in Chapter 7  we give some direct manipulations relating quadratic functions defined using cocycles on manifolds and their boundaries, on manifolds and codimension zero submanifolds, and on a pair of manifolds with a map between them.  All these manipulations are formulated in elementary  terms using the cocycle suspension operation $s$ of Chapter 4, but are somewhat clarified by the simplicial homotopy theory of Chapter 5.  In \S7.3 we also discuss some specific facts about quadratic functions in low dimensions.  Relations between quadratic functions and Pin structures are well-known in dimension 2.  See for example [4] for an elaborate discussion.  In fact, in that dimension, quadratic functions are actually defined and studied rather geometrically on cohomology groups, not  on cocycle groups. Our discussion provides some additional simplicial and homotopy theoretic perspective even in dimension 2. \\

\subsection{Ordered Simplicial Complexes}

Throughout we will work with simplicial sets and with spaces that are {\em ordered} simplicial complexes, that is, the vertices are partially ordered so that the vertices of each simplex are totally ordered.  Maps will be (weakly) order preserving simplicial maps.  Given any simplicial complex $X$, one can always just totally order the vertices.  Given any simplicial map $f \colon Y \to X$ with a (partial) vertex order on $X$, one can easily partially order the vertices of $Y$  so that the map is order preserving.  For example, totally order the vertices in $f^{-1}(v)$, for each vertex $v \in X$, and then arrange  that if $v < w$ in $X$ then the vertices in $f^{-1}(v)$ precede the vertices in $f^{-1}(w)$.  It is an important point that the barycentric subdivision of a simplicial complex has a canonical vertex ordering, where each barycentric vertex is assigned the dimension of its underlying simplex.  The barycentric subdivision of any simplicial map  is always order preserving for the canonical barycentric orderings.  Moreover, if a complex is given any vertex ordering, there is a canonical order preserving simplicial map homotopic to the identity from the barycentric subdivision to the complex.  One maps each barycentric vertex to the maximum vertex of its underlying simplex.\\

Thus, from now on, simplicial complex will always mean ordered simplicial complex and simplicial map will always mean ordered simplicial map.  A simplicial complex determines a simplicial set with $q$-cells named by weakly increasing vertex strings $(v_0v_1...v_q)$ that span a simplex.  Above the dimension of the complex, these are all degenerate.\\

Historically, simplicial complexes and more general cell complexes were viewed as convenient ways to break spaces up into simple pieces, which are glued together following certain rules.  But even defining homology, with anything other than $\Z/2$ coefficients, required more structure on the pieces, for example orientations.  To properly understand the full structure of cohomology, including cup products and cohomology operations, even more is required, for example a diagonal approximation, which can be obtained from a vertex ordering.  Topological and homotopy  invariance was non-trivial, and more or less required properties of barycentric subdivisions.  Still, it seemed to be the case that the simplicial structures with vertex orderings were just a means to an end.  One had to explain why various functors did not depend on these choices.\\

The introduction of the singular complex in topology magically took care of many combinatorial issues, including   topologically invariant definitions of homology groups, cohomology rings, and cohomology operations.  But barycentric subdivisions still played a role, for example in the proof of the excision axiom, and in a sort of a hidden way in the theory of abstract Kan simplicial sets.  If one believes that the study of finer and finer arbitrary triangulations of  spaces is necessary for an understanding of physics, then one is sort of stuck with the view that  combinatorial structures are only a means to an end.  But it seems like a possible alternative is that simplicial structures,  enhanced with canonical vertex orderings related to dimensions of things associated to vertices, is the end.  Then one could regard vertex orderings as a blessing, not a nuisance.  Multivariable cochain operations that exploit the vertex ordering are themselves important, not just their cohomological shadows.\\

We  work with {\em normalized} cochain complexes, consisting of  cochains that vanish on degenerate simplices.  Thus  cochain and cocycle always means normalized cochain and cocycle.  If a simplicial complex $X$ is fixed or understood, we will  write $C^*( F), Z^*(F), H^*(F)$ to indicate  cochains, cocycles, and cohomology of $X$ with coefficients in an abelian group $F$.\\ 

Ordered simplicial structures allow one to define cup products and higher cup$_i$ products in various cochain complexes.  We use the standard formulas of Steenrod for the cup$_i$ products, and we make use of standard properties of these products, especially the coboundary formula.  For integral cochains $X, Y$ of degrees $|X|, |Y|$ that formula is: 
\begin{eqnarray*}
\lefteqn{d(X \cup_i Y) = }\\
 & &(-1)^i\bigl ( dX \cup_i Y + (-1)^{|X|}X \cup_i dY - X \cup_{i-1} Y - (-1)^{i+|X||Y|} Y \cup_{i-1}X \bigr).
\end{eqnarray*}
The $k$-invariants and product formulas and relations between cochain tuples that we encounter in  simple Postnikov towers are expressed in terms of $\cup_i$ operations. \\

For $\Z/2$ cochains $c$, we repeatedly encounter the  cochain definition of the Steenrod  square operations. That definition is, for a cochain $c$ of degree $k$, $$Sq^i c =  c \cup _{k-i} c + c \cup_{k-i+1} dc.$$ If $dc = 0$, this is the standard formula for $Sq^i$ on cocycles.  The cochain formula has the nice property $Sq^i(dc) = dSq^i (c)$.  The shorthand notation $Sq^i c$ doesn't  help with actual computation, but it makes expressions of many of our formulas and proofs of statements about those formulas  more efficient.\\

For example, the coboundary formula leads to the following useful non-linearity relation for $\Z/2$ cochains $c, c'$ of degree $k$.
\begin{equation}
Sq^2(c'+c)   = Sq^2c' +Sq^2c + dc' \cup_{k} dc + d(c' \cup_{k-1} c + dc' \cup_{k} c).
\end{equation}
If $c'$ is a cocycle (1.2) simplifies to
\begin{equation}Sq^2(c' + c) = Sq^2c' + Sq^2c + d(c' \cup_{k-1} c).
\end{equation}

\section {Quadratic Functions}

\subsection {Alternate Definition of Quadratic Functions}

Fix a compact $n$-manifold $(M, \partial M)$, with a simplicial structure.  The boundary can be empty. For later use, we assume the boundary is a full subcomplex and we assume the ordered simplicial structure on $M$ has the property that for each simplex that intersects the boundary, the boundary vertices precede the non-boundary vertices.  This can always be accomplished by passing to one barycentric subdivision and taking the canonical  ordered structure.   We have defined $\Z /4$-valued quadratic functions in the Preview section 1.1.  However, for some purposes it will be more natural to view the values  in $\R / \Z$.\\

By a quadratic  function on $M$ we mean a function  $$Q\colon Z^{n-1}(M, \partial M; \Z/2) \to  \R / \Z$$ that satisfies $$Q(p+q) = Q(p) + Q(q) + (1/2)\int_{[M, \partial M]}\  p \cup_{n-2} q $$ and $$Q(dc) = (1/2)\int_{[M, \partial M]}\ Sq^2c\ = (1/2)  \int_{[M, \partial M]}\  c \cup_{n-4} c + c \cup_{n-3} dc.$$ The $(1/2)$ means the inclusion $\Z/2 \subset \R / \Z$. Note  the pairing $\int_{[M, \partial M]}\ p \cup_{n-2} q$ is symmetric because of the coboundary formula for cocycles $d(p \cup_{n-1} q) = p \cup_{n-2} q + q \cup_{n-2} p$.   Taking $p = q$ and using the fact that $p \cup_{n-2} p = Sq^1p$, we see that $M$ orientable  actually forces $0 = Q(2p) = 2Q(p)$.  In the general Pin case, we see that $0 = Q(2p) = 2Q(p) + (1/2) \int_{[M, \partial M]}\ Sq^1p$, hence $0 = 4Q(p)$ is forced.\\

We already know that the existence of a quadratic function requires $v_2(M) = 0$, which we will assume throughout.  In \S1.1, just after the statement of Claim 1.1(P), we made some other elementary assertions about the definition of quadratic functions, but we postponed the proofs.  We want to check that the formula for $Q(dc)$ depends only on $dc$.  If $x \in Z^{n-2}(M, \partial M; \Z/2)$ is a relative cocycle then from Equation (1.3) $$Sq^2(x + c) = Sq^2x + Sq^2c + d(x \cup_{n-3} c).$$  But $\int_{[M, \partial M]}\ Sq^2x = 0$ since $v_2(M) = 0$ and $\int_{[M, \partial M]}\ d(x \cup_{n-3} c) = 0$ by Stokes Theorem, since $x$ and $c$ vanish on $\partial M$.  Thus $\int_{[M, \partial M]}\ Sq^2(x+c) = \int_{[M, \partial M]}\ Sq^2c$. Similarly, Equation (1.3) implies immediately by Stokes Theorem that the two formulas for $Q(dc + dc')$ given by the two conditions on $Q$ also coincide.\\

Following the statement of Claim 1.1(S) in \S1.1 we gave a construction of quadratic functions, which is now valid since the details of the consistency of the definition of quadratic functions have been established. We want to use that construction to relate quadratic functions on a triangulated $n$-manifold $M$ with a vertex partial ordering to quadratic functions on the first barycentric subdivision $M'$, with its standard vertex partial order.  The key is the observation made in \S1.3 that there is a canonical order preserving simplicial map $b \colon M' \to M$ homotopic to the identity.  The map assigns to each vertex  $v \in M'$ the maximum vertex of the simplex of $M$ containing $v$ as its barycenter.  Any map in the ordered simplicial category induces a map of cochain complexes that commutes with $\cup_i$ products.  In our case, $b^* \colon C^*(M; \Z/2) \to C^*(M'; \Z/2)$ induces a cohomology isomorphism.  Therefore, cocycles $\{p_j\}$ on $M$ defining  a cohomology basis in dimension $n-1$ yield by pull-back a cohomology basis $\{p_j'\}$ on $M'$. Examining the construction of quadratic functions back in \S1.1, we see that setting $Q'(p_j') = Q(p_j) \in \R/ \Z$ establishes the following.\\

RESULT 2.1:  The canonical ordered simplicial map $b \colon M' \to M$ determines a bijective correspondence between quadratic functions on $M'$ and $M$, by means of compositions $$Q = Q' b^* \colon Z^{n-1}(M, \partial M ; \Z/2) \to Z^{n-1}(M', \partial M'; \Z/2) \to \R / \Z.$$
Of course this result also explains why there is a canonical correspondence between quadratic functions on $M$ with two different vertex orders.  They share the same $M'$.  In \S3.4 we will explain using a different method how there is a canonical bijective correspondence between quadratic functions on $M$ with any two simplicial structures.

\subsection {A Construction for Spin Manifolds}

It is obvious that a quadratic function on a disconnected manifold is simply a  sum of quadratic functions on the separate components, so we might as well assume that we begin with a connected Spin manifold $M^n$.  For simplicity, we begin with the closed case, and we assume $M^n$ is a Spin boundary. We will construct a canonical quadratic function $Q$ on $M^n$, following a method of Kapustin, [3].\\

Step 1:  From the Atiyah-Hirzebruch spectral sequence, the reduced $n$-dimensional Spin bordism of the Eilenberg-MacLane space  $K(\Z/ 2, n-1)$ vanishes.  Thus, given $p \in Z^{n-1}(M;\Z/2)$ we can find a Spin manifold $W$ with $\partial W = M$, and a cocycle $\tilde{p} \in Z^{n-1}(W; \Z/2)$ that restricts to $p$ on the boundary.  We will define $$Q(p) = (1/2)\int_{[W, \partial W]}\ Sq^2 \tilde{p} \in \R / \Z.$$

Step 2: We need to prove that $Q(p)$ is well-defined.  First, if $\tilde{p}$ is replaced by $\tilde{p} + \tilde x$, where $\tilde x \in Z^{n-1}(W, M; \Z /2)$, then $$Sq^2(\tilde{p}+\tilde x) = Sq^2\tilde{p} + Sq^2\tilde x + d(\tilde{p} \cup_{n-2} \tilde x).$$  But $W$ is a Spin manifold and $\tilde x$ is a relative cocycle, so $\int_{[W, \partial W]}\ Sq^2\tilde{x} = 0$.  Also, the restriction of $\tilde p \cup_{n-2} \tilde x$ to $\partial W = M$ vanishes.  Thus by Stokes Theorem the integrals over $W$ of $Sq^2 \tilde p$ and $Sq^2(\tilde p + \tilde x)$ coincide.\\

Next, if we replace $W, \tilde p$ by another Spin manifold choice $W',  \tilde  p'$ with boundary $M, p$, then we can form a closed Spin manifold and a cocycle $ \hat{W}, \hat p $, by gluing $W$ and $-W'$ together along $M$.  Since any $n$-simplex of $\hat W$ lies either in $W$ or $W'$, we see that $0 = \int_{\hat W}\ Sq^2 \hat p = \int_{[W, \partial W]}\ Sq^2\tilde p + \int_{[W', \partial W']}\ Sq^2 \tilde p'.$\\

Step 3:  We need to prove that $Q$ is quadratic.  If $c \in C^{n-2}(M; \Z/2)$ then we can lift $c$ to a cochain $\tilde c$ on any Spin manifold $W$ with boundary $M$.  Then $d \tilde c$ lifts $dc$ and $$Q(dc) = (1/2) \int_{[W, \partial W]}\ Sq^2d\tilde c = (1/2)\int_{[W, \partial W]}\ d Sq^2 \tilde c = (1/2)\int_M\ Sq^2c.$$  This verifies the second defining condition for a quadratic function.\\

To verify the quadratic property for a sum $p+q$, we use the fact that the reduced Spin bordism of $K(\Z,2, n-1) \times K(\Z/2, n-1)$ vanishes.  Then we can choose one Spin manifold $W$ with $\partial W = M$, and lifts of both $n-1$ cocycles $p, q$ to cocycles $\tilde p, \tilde q$ on $\tilde W$.  Then from Equation (1.3), $ Sq^2(\tilde p + \tilde q) = Sq^2 \tilde p + Sq^2 \tilde q + d(\tilde p \cup _{n-2} \tilde q)$.  Integrating over $W$ and using Stokes Theorem again gives the desired result $$Q(p+q) = Q(p) + Q(q) + (1/2)\int_M\ p \cup_{n-2}q.$$\\
We need to extend the construction above to Spin manifolds with boundary $(M, \partial M)$, and to closed Spin manifolds $M$ that are not Spin boundaries.  In both cases, we form the closed Spin boundary $DM = M \cup_{\partial M} -M$.  If $\partial M = \emptyset$ then $DM$ is just two copies of $M$ with opposite orientations.  We can use the same vertex order on simplices of $-M$ that are used on the corresponding simplices of $M$. Cocycles $p \in Z^{n-1}(M, \partial M; \Z/2)$ extend by 0 on $-M$ to cocycles $p^+ \in Z^{n-1}(DM; \Z/2)$. We then apply the construction above for closed manifolds to get a quadratic function $Q^+$ on $DM$ and define $Q(p) =  \Q^+(p^+)$.  The steps showing that $Q^+$ is quadratic can be copied to prove that $Q\colon Z^{n-1}(M, \partial M; \Z/2) \to \R / \Z$ is also quadratic.\\

\subsection {A Second Construction for Spin Manifolds}

Given a Spin manifold $(M^n, \partial M^n)$ we will give an alternate construction of the canonical quadratic function $Q\colon Z^{n-1}(M, \partial M; \Z/2) \to \R /\Z$,   at least after possibly subdividing the simplicial structure on $M$. We rely on the Result 2.1 that quadratic functions on $M$ and on any iterated barycentric subdivision $M'$ correspond canonically. So  if we subdivide $M$ and define a quadratic function $Q'$ on $M'$, by Result 2.1 $Q'$ will correspond to exactly one quadratic function $Q$ on the original simplicial structure on $M$. The key now  is that every cohomology class $\bar p \in H^{n-1}(M, \partial M; \Z/2)$ has the form $\bar p = u^*(\bar z)$ for some map $u\colon (M, \partial M) \to (S^{n-1}, v_0)$, where $z$ is a fundamental cocycle on $S^{n-1} = \partial \Delta^n$ that takes value 1 on the face opposite the first vertex, $v_0$, of the simplex $\Delta^n$, and takes value 0 on other faces.  After possibly subdividing $M$, we can assume $u$ is a simplicial map $(M, \partial M) \to (S^{n-1}, v_0)$, with cocycle $p = u^*(z) + dc$ for some $c \in C^{n-2}(M, \partial M)$.  By transversality, $Z =  u^{-1}(pt) \subset M - \partial M$ is a framed, hence Spin, 1-manifold in $M$, where $pt$ is the barycenter of the face  of $\Delta^n$ opposite $v_0$. As such, $Z$ represents an element $[Z] \in \Omega^{spin}_1(pt) = \Z/2$.  We can also directly interpret $[Z]$ as the (relative) Spin bordism class of $u\colon (M, \partial M) \to (S^{n-1}, v_0).$  We then define $Q(u^*(z)) = (1/2)[Z] \in \R / \Z$ and, necessarily, $$Q(p) = (1/2)[Z] + (1/2)\int_{[M, \partial M]}\ Sq^2c + p \cup_{n-2} dc.$$  If we subdivide $M$ enough times, we can assume a finite set of cocycles $p_j$ representing a basis of $H^{n-1}(M, \partial M; \Z/2)$ have the geometric form $u_j^*(z)$. Then the canonical quadratic function is determined by $Q(p_j) = (1/2)[Z_j] \in \R / \Z$, along with the two defining conditions for a quadratic function.  \\

In ([1], \S6.3) we gave an argument that the two formulas for $Q(p)$ coincide for $n = 3$.  Here is an easier argument for $n > 3$, which we are reproducing from ([2], \S10.3).  A little computation shows we may as well assume $p = u^*(z)$ for a simplicial map $u \colon M \to S^{n-1}$, where $M$ is a closed Spin boundary.  The reduced $n$-dimensional Spin bordism of $\Sigma^{n-3}\mathbb{C}P(2)$ vanishes, so we can choose a simplicial structure on $\Sigma^{n-3}\mathbb{C}P(2)$ which contains $S^{n-1}$ as a subcomplex, and, after possibly subdividing, extend $u$ to a simplicial map $$\tilde u \colon (W, M) \to (\Sigma^{n-3}\mathbb{C}P(2), S^{n-1}),$$ where $W$ is a Spin manifold with $\partial W = M$.  Then $\tilde u^*(\tilde z) = \tilde p \in Z^{n-1}(W; \Z/2)$ restricts to $p$ on $M$, where $\tilde z$ is a cocycle representing the generator of $$H^{n-1}(\Sigma^{n-3} \mathbb{C}P^2; \Z/2).$$  The cocycle $Sq^2 \tilde z = \tilde z \cup_{n-3} \tilde z$ is a relative cocycle for $(\Sigma^{n-3}\mathbb{C}P(2), S^{n-1})$ that represents a generator of $H^{n+1}(\Sigma^{n-3}\mathbb{C}P(2), S^{n-1}; \Z/2)  = \Z/2$.  This is the key. On the one hand, the integral of this cocycle pulled back to $W$ in the top dimension coincides with $\int_{[W, \partial W]}\ Sq^2 \tilde p$.  On the other hand, this integral is exactly the obstruction to deforming $\tilde u \colon (W, M) \to (\Sigma^{n-3}\mathbb{C}P(2), S^{n-1})$ rel $M$ to a map $W \to S^{n-1}$, and this obstruction is the Spin bordism class of $u\colon M \to S^{n-1}$.  Thus the two formulas for $Q(p)$ coincide.\\

For the most part in this section, we have in mind that the cocycles $p$ represent non-trivial cohomology classes.  But the results apply to any cocycle.  For example, if we  look at the link of an $n-2$ simplex $\sigma$ in $M$, we see a small framed circle $Z$ dual to $d\sigma^*$ bounding a framed disk, hence the link inherits the trivial Spin structure from $M$.  This `explains' the formula $Q(d\sigma^*) = 0$ that we asserted back in \S1.1 was equivalent to the defining formula for $Q(dc)$ for all coboundaries.  One obtains all $dc$ by adding various $d\sigma^*$.\\
 
 \subsection {Change of Spin Structure on a Manifold}
 
 Suppose $a \in Z^1(M; \Z/2)$ and suppose $M_a$ denotes the new Spin structure on $M$ obtained by the action of $H^1(M; \Z/2)$ on Spin structures by classical mechanisms.  We want to prove the new canonical quadratic function $Q_a$ on $Z^{n-1}(M, \partial M; \Z/2)$ is given by $$Q_a(p) = Q(p)+ (1/2)\langle ap, [M, \partial M] \rangle.$$
To study $Q_a - Q$, without loss of generality  we can use Result 2.1 and subdivide and assume $p = u^*z$, where $u\colon (M, \partial M) \to (\partial \Delta^n, v_0)$ is simplicial and $z$ is a fundamental cocycle on the sphere non-zero only on the face opposite $v_0$.   Then $p \in Z^{n-1}(M, \partial M; \Z/2)$ is dual to a framed 1-manifold $Z = u^{-1}(pt) \subset M$, where $pt$ is the barycenter of the face opposite $v_0$.    The two Spin structures $M$ and $M_a$ induce Spin structures $[Z]$ and $[Z_a]$. From our discussion of the quadratic functions defined on $Z^{n-1}(M, \partial M; \Z/2)$ above, $Q_a(p) - Q(p)$ will be the difference $[Z_a] - [Z]$.    A Spin structure on $M$ can be interpreted as a choice of trivialization of the bundle $\tau_M + 1$ on a neighborhood of a dual 2-skeleton of $M$. Here $\tau_M$ is the tangent bundle of $M$.  Changing the Spin structure on $M$ by the cohomology class of $a$ changes the  trivialization of the restriction of $\tau_M + 1$ to $Z$ by the homotopy class of a map $$Z \to M - M^{(n-3)} \xrightarrow{a} \R P^n \to SO(n+1).$$  As a homology class in $\Z/2 = H_1(\R P^n ; \Z/2)$, this  class coincides with the number $\langle a, [Z]\rangle = \langle ap, [M, \partial M] \rangle$.  This proves $$Q_a(p) = Q(p) + (1/2)\langle ap, [M, \partial M]\rangle.$$
Thus, as the Spin structure on $M$ varies, we see a bijective correspondence between Spin structures on $M$ and  quadratic functions $Q$.  We have proved the following result, which is just a restatement of Claim 1.1(S).\\

\noindent RESULT: Fix a compact oriented simplicial $n$-manifold $M$. Equivalence classes of Spin structures on $M^n$ are in canonical bijective correspondence with functions $$Q\colon Z^{n-1}(M, \partial M; \Z/2) \to  \R / \Z$$ that satisfy $2Q(p) = 0$, along with $$Q(p+q) = Q(p) + Q(q) + (1/2)\int_{[M, \partial M]}\  p \cup_{n-2} q $$ and $$Q(dc) = (1/2)\int_{[M, \partial M]}\ Sq^2c\ =  (1/2)\int_{[M, \partial M]}\  c \cup_{n-4} c + c \cup_{n-3} dc.$$

\subsection {A Construction for Pin Manifolds}

Now let $(M^n, \partial M^n)$ be a (compact) Pin manifold.  We want to at least indicate the definition of a canonical quadratic function $Q \colon Z^{n-1}(M, \partial M; \Z/2) \to \R / \Z$.  Consider $p \in  Z^{n-1}(M, \partial M; \Z/2)$.  Elementary obstruction theory implies that there is a map $u\colon (M, \partial M) \to (\Sigma^{n-2}\R P^2, v)$ with $u^*(\bar z) = \bar p$, where $v$ is a base point and $z$ is a cocycle generating $H^{n-1}(\Sigma^{n-2}\R P^2; \Z/2)$. After possible subdivision, we can assume $u$ is simplicial and $u^*(z) = p + dc$, for some cochain $c \in C^{n-2}(M, \partial M; \Z/2)$.\\

Let $E^2 \subset M - \partial M$ be a transverse inverse image under $u$ of the framed submanifold $\R P^2 \subset  \Sigma^{n-2}\R P^2$.  Then $E^2$ is a closed Pin surface, which is a Pin boundary\footnote{If $\Sigma X = C^+X \cup_X C^-X$ is expressed as a union of two cones, and if $f \colon (M, \partial M) \to (\Sigma X, v_-)$ is a map transverse to $X$ where $v_-$ is the lower cone vertex, then $N =  f^{-1}(X)$ is the boundary of $f^{-1}(C^+X) \subset M$.  Iterate to understand maps to higher suspensions.} with a map $E^2 \to \R P^2$.  It is fairly well-known that Pin structures on a surface $E^2$ correspond canonically to quadratic functions $q \colon H^1(E^2; \Z/2) \to \Z / 4$ refining the cup product pairing, [4].  Also, the reduced Pin bordism of $\R P^2$ is isomorphic to $\Z/4$, (this will be proved in \S6), and the Pin bordism class of a map $u \colon E^2 \to \R P^2$ is the value $q(u^*(x)) \in \Z/4$, where $x \in H^1(\R P^2; \Z/2)$ is the generator.\\

Now, back to our map $u\colon (M, \partial M) \to (\Sigma^{n-2}\R P^2, v)$.  The canonical quadratic function on $M$ will be defined by $Q(u^*(z)) = (1/4)[E^2 \to \R P^2] \in (1/4)\widetilde{\Omega}^{Pin}_2(\R P^2) \subset \R / \Z$.  Then, with $p = u^*(z) + dc$, we (necessarily) define $$Q(p) = (1/4)[E^2 \to \R P^2] + (1/2)\int_{[M, \partial M]}\ Sq^2 c + p \cup_{n-2} dc \in \R /\Z.$$
We do not have a direct proof that $Q$ is a quadratic function in the Pin case.  The difference between the Spin and Pin cases is that we do not have an analogue of the Kapustin method for defining $Q$.  The reduced Pin bordism of $K(\Z/2, n-1)$ does not vanish. Even if an $n-1$ cocycle $p$ on $M^n$ does extend to a cocycle on a Pin manifold $W$ bounding $M$, it is unclear how to use $W$ to define a value $Q(p) \in \R / \Z$, which must have order 4 if $Sq^1p \not= 0$.\\

We will therefore develop a more indirect homotopy theoretic method that works in both the Spin and Pin cases.  The idea is to translate quadratic functions on $M^n$ to linear characters on groups $G_n(M^n, \partial M^n)$ that are Pontrjagin duals of quotients of Spin or Pin bordism.  Then, by duality,  the identity bordism element $Id \colon (M^n, \partial M^n) \to (M^n, \partial M^n)$ will automatically define a linear function $G_n(M^n, \partial M^n)) \to \R / \Z$, which is  equivalent to a quadratic function on $M^n$.  We can prove in the Pin case that this quadratic function is indeed given by the formula for $Q(p)$ just above.

\section {The Groups $\Gsp$ and $\Gp$}

\subsection{Definition of $\Gsp$}

Fix a simplicial complex or simplicial set $X$.  Begin with pairs $$\C^{spin}_n(X) = \{(w, p) \in C^n(\R / \Z) \times C^{n-1}(Z/2) \vert \ dp = 0,\ dw = (1/2)Sq^2p\},$$ where $(1/2)$ means the coefficient morphism $\Z/2 \to \R / \Z$.  Define a product by $$(w, p)(v, q) = (w+v + (1/2)p \cup_{n-2}q,\ p + q).$$  The product is  associative,  has an identity element $(0, 0)$, and the inverse of $(w, p)$ is  $(-w + (1/2)Sq^1p,\ p)$.  Thus we have a group. The product is not commutative, since in general $p \cup_{n-2}q \not= q \cup_{n-2}p$.  However, we can understand all commutators from the formula $$(w,p)(v,q) =(v,q)(w,p) (d((1/2)p \cup_{n-1} q),\ 0).$$
We will divide this group of pairs $(w, p)$ by the subgroup consisting of elements $(df + (1/2)Sq^2c, dc)$. This is a subgroup because of the relation (1.2) from \S1.3, $$  Sq^2(c'+c)  = Sq^2c' +Sq^2c + dc' \cup_{n-2} dc  + d(c' \cup_{n-3} c + dc' \cup_{n-2} c).$$  Since the subgroup of elements $(df, 0)$ already contains all commutators, the quotient group $$\Gsp = \C^{spin}_n(X) / \{ (df + (1/2)Sq^2c, \ dc)\}$$ is an abelian group.  Since $Sq^1p$ is the reduction of an integral torsion class, the $\R/\Z$ cocycle  $(1/2)Sq^1p$ is a coboundary.  Thus after dividing by relations, the inverse of $(w, p) \in \Gsp$ is given simply by $(-w, p)$ \\

Directly from the definitions it is easy to see that the maps on the cochain level  $w \mapsto (w, 0)$ and $(w,p) \mapsto p$ induce a short exact sequence
\begin{equation}
0 \to QH^n(X; \R /\Z) \to \Gsp \to SH^{n-1}(X; \Z/2) \to 0,
\end{equation}
where $$QH^n(X; \R /\Z) = H^n(X; \R / \Z)\ /\ Image(H^{n-2}(X; \Z/2) \xrightarrow{(1/2)Sq^2} H^n(X; \R / \Z))$$ and  $$SH^{n-1}(X; \Z/2) = Kernel( H^{n-1}(X; \Z/2) \xrightarrow{(1/2)Sq^2} H^{n+1}(X ; \R / \Z)).$$

\subsection{Definition of $\Gp$}

Next we turn to the Pin version.  Begin with pairs $$\C^{pin}_n(X) = \{(w, p) \in C^n(\Z/2) \times C^{n-1}(Z/2) \vert \ dp = 0,\ dw = Sq^2p\}.$$
Define a product by $$(w, p)(v, q) = (w+v + p \cup_{n-2}q,\  p + q).$$
The product is  associative,  has an identity element $(0, 0)$, and the inverse of $(w, p)$ is  $(-w + Sq^1p,\ p)$.  Thus we have a group. The product is not commutative, since in general $p \cup_{n-2}q \not= q \cup_{n-2}p$.  However, we can understand all commutators from the formula $$(w,p)(v,q) = (d(p \cup_{n-1} q), 0)(v,q)(w,p).$$
We will divide this group of pairs $(w, p)$ by the subgroup consisting of elements $(df + Sq^2c, dc)$. Just as in the Spin case, this is a subgroup because of the relation (1.2) from \S1.3.  Since the subgroup of elements $(df, 0)$ already contains all commutators, the quotient group $$\Gp = \C^{pin}_n(X) / \{ (df + Sq^2c, \ dc)\}$$ is an abelian group.\\

Directly from the definitions it is easy to see, as in the Spin case,  that the maps on the cochain level  $w \mapsto (w, 0)$ and $(w,p) \mapsto p$ induce a short exact sequence
\begin{equation}
0 \to QH^n(X; \Z/2) \to \Gp \to SH^{n-1}(X; \Z/2) \to 0,
\end{equation}
where $$QH^n(X; \Z/2) = H^n(X; \Z/2)\ /\ Image(H^{n-2}(X; \Z/2) \xrightarrow{Sq^2} H^n(X; \Z / 2))$$ and  $$SH^{n-1}(X; \Z/2) = Kernel( H^{n-1}(X; \Z/2) \xrightarrow{Sq^2} H^{n+1}(X ; \Z/2)).$$
In both the Spin and Pin versions, the cochain and cocycle definitions of the groups $G$ make sense for pairs $(X, Y)$ of simplicial sets, as well as for single spaces $X$.  It is obvious from the product formulas that the map on cochain representatives $(w, p) \mapsto ((1/2)w, p)$ induces a group homomorphism $G^{pin}_n(X, Y) \to G^{spin}_n(X, Y)$.

\subsection{Homotopy Invariance of $\Gsp$ and $\Gp$}

The groups $\Gsp$ and $\Gp$  obviously define contravariant  functors on the (ordered) simplicial category.  Moreover, the exact sequences (3.1) and (3.2) are functorial  invariants.  Suppose $f\colon X \to Y$ is a simplicial map inducing isomorphisms in cohomology with $\Z/2$ and $\R / \Z$ coefficients.  Since the induced cohomology maps commute with $Sq^2$, we see from the exact sequences (3.1) and (3.2)  that the induced maps $f^*\colon G^{spin}_n(Y) \to G^{spin}_n(X)$ and $f^*\colon G^{pin}_n(Y) \to G^{pin}_n(X)$ are isomorphisms.  The same statements apply to maps between pairs of spaces. \\

In particular, the remarks in the above paragraph apply to homotopy equivalences.  For example, if $I \times X$ is given a simplicial structure, restricting to simplicial structures $X_0$ and $X_1$ on the ends $\{0\} \times X$ and $\{1\} \times X$ respectively, then the inclusions $X_i \to I \times X$, $i = 0,1$,  induce isomorphisms between all $G_n$ groups.  This is one way to express the result that the $G_n$ groups do not depend on the simplicial structure on a space.  But we want to free these isomorphisms from  choices of simplicial structures on $I \times X$.\\

We can make the invariance more precise by exploiting the singular complex $S(X)$.  Given any simplicial structure on a space $X$, there is a canonical simplicial map $X \to S(X)$, which induces an isomorphism in cohomology groups, and hence also in the $G_n$ groups.  Now the projection of topological spaces $I \times X \to X$ is not necessarily a simplicial map, but it does induce a simplicial map $S(I \times X) \to S(X)$. We can compose this last map with the simplicial inclusions $X_i \to I \times X \to S(I \times X)$, $i = 0,1$.  The conclusion is that not only are the groups $G_n(X_0)$ and $G_n(X_1)$ isomorphic, but they are both canonically identified with $G_n(S(X))$, in both the Spin and Pin cases.  \\ 

As an additional consequence of bringing in the singular complexes of $X$ and $I \times X$, we can conclude that if $f_0 \colon X_0 \to Y_0$ and $f_1 \colon X_1 \to Y_1$ are topologically homotopic simplicial maps between  simplicial structures on $X$  and simplicial structures on $Y$, then the induced maps $f_i ^* \colon G_n(Y_i) \to G_n(X_i)$ for $i = 0,1$ coincide, when these $G_n$ groups are identified with $G_n$ groups of the singular complexes of $X$ and $Y$ respectively.\\

The paragraphs above establish in a strong way that the groups $\Gsp$ and $\Gp$ are homotopy functors.  This will become  more transparent in  Chapter 6, where we will prove that the groups can be described as homotopy classes of maps from $X$ to spaces $E^{spin}_n$ and $E^{pin}_n$, respectively.

\subsection{Quadratic Functions as Linear Functions $G_n(M^n) \to \R / \Z$}

We now take $X$ to be a compact simplicial $n$-manifold $M^n$, possibly with boundary.  For dimensional reasons, $Sq^2 p = 0$ and $dw = 0$, so in the $G_n$ groups we are beginning with all pairs $(w, p)$ of relative cocycles of the relevant degrees in both the Spin and Pin cases.  The class $p$ is always a $\Z/2$ cocycle and $w$ is an $\R / \Z$ cocycle in the Spin case, and a $\Z/2$ cocycle in the Pin case.  Moreover, we have the factorization $(w, p) = (w, 0)(0, p)$ in the $G_n$ groups.\\

We will also assume $v_2(M) = 0$, and we will assume, without really losing any generality, that $M$ is connected.  The short exact sequences (3.2) and (3.1) simplify to
$$0 \to H^n(M, \partial M; \Z/2) \to G^{pin}_n(M, \partial M) \to H^{n-1}(M, \partial M; \Z/2) \to 0$$ in the Pin case, and to $$0 \to H^n(M, \partial M; \R / \Z) \to G^{spin}_n(M, \partial M) \to H^{n-1}(M, \partial M; \Z/2) \to 0$$ in the Spin case. The group homomorphism $G^{pin}_n(M, \partial M) \to G^{spin}_n(M, \partial M)$, defined by $(w, p) \mapsto ((1/2)w, p)$,  induces a map from the first short exact sequence to the second.  We have $H^n(M, \partial M; \Z/2) = \Z /2$ and we  have $H^n(M, \partial M; \R/ Z) = 0$ unless $M$ is orientable, in which case it is isomorphic to $\R / \Z$.  In that case there are two (continuous) isomorphisms, defined by integrating over the relative fundamental class, with a choice of orientation.\\  

First consider the Spin case.  From the product formula $(w, p)(v, q) = (w+v + (1/2)p \cup_{n-2} q,\ p+q)$ and the factorization $(w, p) = (w, 0)(0, p)$, the following claim is obvious.\\

\noindent CLAIM 3.1:  If $M$ is oriented there is a canonical bijective correspondence between quadratic functions $Q \colon Z^{n-1}(M, \partial M; \Z2) \to \R / \Z$ and linear functions 
$L\colon G^{spin}_n(M, \partial M) \to \R / \Z$ satisfying $L(w, 0) =  \int_{[M, \partial M]}\ w$.  The correspondence is given by  $Q \leftrightarrow L_Q$, where $L_Q(w, p) = Q(p) +  \int_{[M, \partial M]}\ w $. \\

For an arbitrary manifold,  the Pin case is equally obvious from the product formula $(w,p)(v, q) = (w+v + p \cup_{n-2}q,\ p+q)$ and the factorization $(w, p) = (w, 0)(0, p).$\\

\noindent CLAIM 3.2: For any manifold, there is a canonical bijective correspondence between quadratic functions $Q \colon Z^{n-1}(M, \partial M; \Z2) \to \R / \Z$ and linear functions $L\colon G^{pin}_n(M, \partial M) \to \R / \Z$ satisfying $L(w, 0) =  (1/2)\int_{[M, \partial M]}\ w$.  The correspondence is given by $Q \leftrightarrow L_Q$, where $L_Q(w, p) = Q(p) +  (1/2)\int_{[M, \partial M]}\ w$.  \\

We remind again that in the Spin case the quadratic functions have values in $\Z/2 \subset \R / \Z$, and in the general Pin case the quadratic functions have values in $\Z/4 \subset \R / \Z$.\\

In the previous section we have explained how as the ordered simplicial structure on $M$ changes, all the groups $G_n(M, \partial M)$ are canonically identified with one another, in both the Spin and Pin cases.  Therefore when we define  `canonical' quadratic functions, via certain linear functions on $G_n(M, \partial M)$, it will be implicit that we are defining the same quadratic function for all simplicial structures, when the $G_n$ groups are identified.\\

In the next three chapters, we will explain how the groups $\Gsp$ and $\Gp$ are related to Pontrjagin duals of Spin and Pin bordism groups of $X$.  Then, regarding an identity map $Id \colon (M, \partial M) \to (M, \partial M)$ as a  Spin or Pin bordism element, we get by duality canonical linear functions $G_n(M, \partial M) \to \R / \Z$ from a Spin or a Pin structure on a manifold $M^n$.  Thus, from the Claims 3.1 and 3.2, we get quadratic functions from Spin or Pin structures.  In the Spin case, all this will just amount to a reformulation of the proof of main theorem Claim 1.1(S) given in Chapter 2.  But in the Pin case, we obtain in this rather roundabout way the only proof of Claim 1.1(P) that we know.\\

\section {Suspension of Cochains} 

\subsection{Suspension of Cochains}

We regard the suspension $\Sigma X$ of a space to be the obvious  union of two cones $C^{+}X$ and $C^{-}X$.  Given a triangulation of $X$ with vertex order,  we label the new upper cone vertex $+ \infty$ and the new lower cone vertex $-\infty$.  That is, the cone vertices only occur as the last or first vertex of simplices in the suspension.  But we won't use the lower cone vertex. Given a cochain $c \in C^n(X)$, with any coefficients, we define $s(c) \in C^{n+1}(\Sigma X)$ as follows.  On any simplex in the lower cone $C^{-}X$, the value of $s(c)$ will be 0.  On a simplex of form $(012 ... n \infty)$ in the upper cone $C^{+}X$, the value will be $s(c)(01 ... n \infty) = c(01 ... n)$.\\

The `cone vertex last' convention results in the easily proved formula $sd = ds$.  That is, $s$ is a map of cochain complexes $$s \colon C^*(X) \to C^{*+1}(\Sigma X, C^- X).$$   The cochain suspension map  $s$ induces suspension isomorphisms on cohomology with any coefficients $\tilde H^*(X) \simeq H^{*+1}(\Sigma X, C^-X) \simeq H^{*+1}(\Sigma X)$.\\

\subsection{Commutativity of Suspension and $\cup_i$ Products}

The suspension $s$ has some very nice properties relating $\cup_i$ products in $X$ and $\Sigma X$. First, we point out that with the given ordered triangulation of $\Sigma X$, all ordinary cup products $sx \cup_0 sy$ are 0.  The reason is, an ordered simplex can have at most one vertex $+\infty$, so a proper `first face' will always lie in $C^{-}X$.  On the other hand, the following remarkable formula holds for all $i \geq 0$:  
\begin{equation}
s (x \cup_i y)= (-1)^{deg(x)+i+1}sx \cup_{i+1} sy.
\end{equation}
We believe this is an important formula.  It is not easy to prove.   Obviously it implies  that Steenrod square operations commute with suspension, not just on cohomology and cocycles, but actually on all cochains.  The cochain $\cup_i$ operations generalize to other multi-variable cochain operations and there should be useful extensions of this suspension formula to these other operations.\footnote{In hindsight, we believe using the lower cone $C^-(X)$ to define a suspension $s$ with $ds = -sd$ is a more natural choice. But then in order to get the cleanest formulas relating $s$ and cup$_i$ products it is necessary to use alternatives to the historical definitions of cup$_i$ products, including ordinary cup product.  But it is hard to overturn historical conventions!}\\

In the next Chapter we will interpret the cochain suspension map $s$ not just as an operation on cochains, but as providing actual simplicial maps between certain  Postnikov towers.  For this, we need the commutativity formula (4.1).\\

\section{Simplicial Models for  Postnikov Towers} 

\subsection{Basic Models}

In this section, we briefly recall some classical algebraic topology concerning the simplicial theory of Postnikov towers.  By introducing this rather generally here, the special cases that we need below for our study of quadratic functions and Pin and Spin structures  are easily seen to fit into a general framework.\\

Suppose $A_1$ is an abelian group.  A simplicial set model for an Eilenberg-MacLane space $E_1 =  K(A_1, n_1)  $ has as the $q$-simplices the set of $A_1$-valued normalized $n_1$-cocycles on the standard $q$-simplex $\Delta^q$, with the standard face and degeneracy operations. Note $K(A_1, n_1)$ has a tautological fundamental cocycle  $ a \in Z^{n_1}(K(A_1, n_1); A_1)$.  Also, $K(A_1, n_1)$ is a simplicial abelian group. If $X$ is a simplicial complex, or if $X$ is a simplicial set, then a simplicial map $X \to E_1$ is exactly an $A_1$-valued $n_1$-cocycle, say $a$,   on $X$.  Maps $a_0$ and $a_1$ are homotopic if there is a cocycle $\hat a$ on $I \times X$  that restricts to $a_i$ on $\{i\}\times X$, for $i = 0, 1.$ The  null-homotopic maps are those that extend to a simplicial map $CX \to E_1$, where $CX$ is the cone on $X$. This is equivalent to saying $a = dp$ for some $n_1-1$ cochain $p$ on $X$. Thus homotopy classes of simplicial maps $[X, E_1] = H^{n_1}(X; A_1)$. \\

Next, we discuss 2-stage Postnikov towers $$E_2 = K(A_1, n_1) \ltimes_{k(a)} K(A_2, n_2).$$ The notation is meant to indicate a principal fibration over the first Eilenberg-MacLane space with fiber the second Eilenberg-MacLane space.  The term  $k(a)$ is a cocycle representing the cohomology $k$-invariant of the fibration  $$ \bar{k}( a) \in H^{n_2+1}(K(A_1, n_1); A_2).$$ In practice, $k(a)$ is a natural cochain level version of a cohomology operation.  The q-cells of $E_2$ are given by pairs $(p, a)$, where $ a$ is an $A_1$-cocycle on $\Delta^q$ and $p$ is an $A_2$-cochain on $\Delta^q$ with $dp = k(a)$.  Then a simplicial map $X \to E_2$ is given by a pair of cochains $(p, a)$ on $X$, with appropriate coefficients and dimensions, with $da = 0$ and $dp = k(a)$. Maps $(p_0, a_0)$ and $(p_1, a_1)$ are homotopic if there is a pair of cochains $(\hat{p}, \hat{a})$ on $I \times X$ with $d\hat{a} = 0$ and $d\hat{p} = k(\hat{a})$ that restricts to $(p_i, a_i)$ on $\{i\}\times X$, for $i = 0, 1.$ \\

We can continue and define 3-stage Postnikov towers 
$$ E_3 = K(A_1, n_1) \ltimes_{k(a)}  K(A_2, n_2) \ltimes_{k(p, a)} K(A_3, n_3).$$  Thus $E_3$ is a principal  fibration over $E_2$, with an Eilenberg-MacLane space fiber and $k$-invariant determined  by a natural cocycle $k(p, a) \in Z^{n_3+1}(E_2; A_3)$.  More precisely, $k(p,a)$ is a natural cochain that is a cocycle if $da = 0$ and $dp = k(a)$. Simplicial maps $X \to E_3$ are given by triples of cochains $(w, p, a)$ on $X$, with appropriate coefficients and dimensions, with $da = 0$, $dp = k(a)$, and $dw = k(p, a)$. Homotopic triples are defined just as in the 2-stage case.\\

Cohomology groups of Eilenberg-MacLane spaces are known, and we also have explicit cocycle formulas representing cohomology classes.  Thus, simplicial descriptions of 2-stage Postnikov towers $E_2$ are quite explicit.  But the catch is at the next stage. The cohomology groups of $E_2$ are generally obscure.  Even if the cohomology is known in some sense, cocycle formulas for the cohomology classes may be intractable.\\

\subsection{$H$-spaces}

A Postnikov tower represents an abelian group valued homotopy functor if it is a simplicial $H$-space, that is, there is a  homotopy associative and commutative simplicial product $E \times E \to E$, with a simplicial  homotopy inverse. In the 3-stage case, this  means a product of triples $(w, p, a)(v, q, b) = (u, r, c)$, with appropriate properties. In the case of an H-space, the null-homotopic triples determine all relations between triples because of the group structure.  A 3-stage simplicial Postnikov tower $H$-space representing the Pontrjagin dual of reduced 4-dimensional Spin bordism was studied extensively in our paper [2].\\

In this paper, we only work with  2-stage towers $E = K(A_1, n_1) \ltimes_{k(a)} K(A_2, n_2)$, with $A_1 = \Z /2$ and $A_2 = \R / \Z\ \text{or}\ \Z/2$.  A simplicial product  $E \times E \to E$ will be given by a formula $(p,a)(q,b) = (p + q + u(a, b), a+b)$, where $du = k(a+b) - k(a) - k(b)$.  In our case, $u(a, b)$ will be bilinear, hence the product is associative on the level of maps, not just homotopy classes of maps.  However, $u(a, b)$ is not symmetric, so the simplicial product itself is not commutative.  Instead, we will have $u(a, b) - u(b, a) = dc(a, b)$.  Then $(p, a)(q, b) = (q, b)(p, a)( dc(a, b), 0),$ which implies homotopy commutativity since $(dc, 0)$ is null-homotopic.\\

\subsection{Loop Spaces}

Given a Postnikov tower $E$, its loop space should also be a Postnikov tower, with the degrees of the Eilenberg-MacLane spaces shifted down by one.  This  is neatly explained using the cochain suspension operation $s$ of \S4.  Specifically, among the $q+1$ cells of $K(A, n+1)$, which, recall,  are the $A$-valued $(n+1)$-cocycles on $\Delta^{q+1} = C^+(\Delta^q)$, there are the cocycles $sz$, where $z $ is an $n$-cocycle on $\Delta^q$.  These cocycles $sz$ vanish on the $(n+1)$-simplices of $\Delta^q \subset C^+(\Delta^q)$.  Thus $s$ can be viewed as a map $S \colon C^+(K(A, n)) \to K(A, n+1)$, which maps all cells in the base $K(A, n)$ to the basepoint, the cocycle that is identically zero.\\

The key idea is that the expectation $K(A, n) \sim  \Omega K(A, n+1)$ should be realized as the topological adjoint of a map $S\colon C^+( K(A, n))/ K(A,n)  \to K(A, n+1)$.   But this is exactly what the cochain operation $s$ of \S4 does, as explained in the previous paragraph.  Moreover, the map $S$ is a simplicial map, hence induces a map of cochain complexes $S^*$ from $m+1$-cochains on $K(A, n+1)$ to $m+1$-cochains on $C^+K(A, n)$ that vanish on $K(A, n)$. This last is isomorphic, under $s$, to the $m$-cochains on $K(A, n)$.   To help follow this discussion, denote the fundamental cocycle of $K(A, n+1)$ by $\alpha$ and the fundamental class of $K(A, n)$ by $a$.  Then $S^*(\alpha) = sa$.\\

Now we can associate to a two stage Postnikov system 
$$E = K(A_1, n_1) \ltimes_{k(\alpha)} K(A_2, n_2),$$ another two stage Postnikov system
$$\tilde E = K(A_1, n_1 - 1) \ltimes_{k(a)} K(A_2, n_2 -1),$$  where the cohomology operation $k(a)$ is defined by $$sk(a)  = S^*(k(\alpha)) = kS^*(\alpha) = k(sa).$$  Moreover, there is a canonical simplicial map  $S\colon C^+ \tilde E \to E$ that maps the cone on a  $q$-cell $(p,a)$ of $\tilde E$ to the $(q+1)$-cell $(sp, sa)$ of $E$, and maps cells of $\tilde E$ to the basepoint. The main point here is that $dp = k(a)$ implies $d(sp) = s(dp) = sk(a) = k(sa)$.  A topological adjoint of this map $S$ defines a homotopy equivalence $\tilde E \to \Omega{E}$.\\

Next, suppose $E$ is an $H$-space, with a simplicial product.  Then $\tilde E$ inherits two homotopy products, one from the loop functor applied to the product on $E$ and the other because $\tilde E$ is a loop space.  A standard argument proves these two $H$-space structures on $\tilde E$ are equivalent up to homotopy. \\

Suppose the simplicial product on $E$ is given by $$(\rho, \alpha)(\phi, \beta)= (\rho + \phi + \mu(\rho, \phi), \alpha + \beta).$$  Then the simplicial product on $\tilde E$ can be given by $$(p, a)(q, b) = (p+q+ u(a, b), a+b),$$ where the cochain operation $u(a, b)$ is defined by $su(a, b) = \mu(sa, sb).$  The map $S \colon C^+ \tilde E / \tilde E \to E$  can be interpreted as a functorial correspondence from simplicial maps $\{X, \tilde E\}$ to simplicial maps $\{C^+X/ X, E\}$. At the cochain level, this map is defined by $s(p, a) = (sp, sa)$, and at the cochain level the map commutes with products.\\

These ideas suggest the following perhaps useful definition of a simplicial 2-stage Postnikov spectrum.  This will mean a sequence of simplicial  {H}-space Postnikov towers $E_j = K(A_1, n_1 + j) \ltimes_ {k_j(a_j)} K(A_2, n_2 + j), $ together with simplicial maps $S \colon C^+ E_j / E_j \to E_{j+1}$ that preserve the simplicial products in the sense above,  and whose topological adjoints give homotopy equivalences $E_j \to \Omega E_{j+1}$.  Moreover, on cells one should have $SC(p,a) = (sp, sa)$, where $C$ means cone on a cell of $E_j$.  The definition, including the key role played by the cochain suspension map $s$, should extend to 3-stage Postnikov spectra and beyond.\\

\section{The Classifying Spaces for $G^{spin}_n(X)$ and $G^{pin}_n(X)$ }

\subsection{A Simplicial 2-Stage Postnikov Spectrum $E^{spin}$}

Consider the two stage Postnikov tower
$$ E^{spin}_n = K(\Z/2, n-1) \ltimes_{(1/2)Sq^2p} K(\R / \Z, n).$$ 
Here, $p$ is the fundamental $\Z/2$ cocycle of degree $n-1$ and $Sq^2p = p \cup_{n-3} p$ is the standard cocycle representative of the cohomology operation $Sq^2$ in this dimension. The $(1/2)$ means the coefficient morphism $\Z/2 \to \R / \Z$.  We will define a homotopy commutative $H$-space structure on $E^{spin}_n$ and it will be immediate that there is an isomorphism $ G^{spin}_n(X) = [X, E^{spin}_n]$, where  $G^{spin}_n(X) $ is the group constructed in \S3.1. The space $E^{spin}_n$ represents the Pontrjagin dual  of a quotient of reduced $n$-dimensional Spin bordism.  Specifically, there is a natural isomorphism $$ G^{spin}_n(X) = [X, E^{spin}_n]  \simeq {\rm Hom}(\widetilde{\Omega}_n^{spin}(X) / Image\ \widetilde{\Omega}_n^{spin}(X^{(n-2)}),\ \R /\Z) $$ where $X^{(n-2)}$ denotes the $n-2$ skeleton of $X$.  We will explain this isomorphism in the next section.\\

Elements of $[X, E^{spin}_n]$ are represented by pairs $(w, p)$ with $dp = 0$ and $dw = (1/2)p \cup_{n-3} p$.  The product will be given by
$$(w,p)(v,q) = (w+v + (1/2)p \cup_{n-2} q,\  p+q).$$
We have here a special case of the simplicial $H$-space discussion of \S5.2. Note $$d(p \cup_{n-2} q)  = p \cup_{n-3}q + q \cup_{n-3} p  = Sq^2(p+q) - Sq^2(p) - Sq^2(q).$$ 
The product is strictly associative, but not strictly commutative.\\

The null-homotopic pairs turn out to be $(df + (1/2) Sq^2c,\ dc)$.  We are then exactly in the situation describing the group $G^{spin}_n(X) $ studied in \S3. In particular, the product is homotopy commutative and the inverse of $(w, p)$ is $(-w, p)$. 
Thus pairs modulo all null-homotopic pairs form an abelian group $[X, E^{spin}_n]$, and this group is exactly the group $G^{spin}_n(X) $ studied in \S3.1\\

As mentioned, the arguments above can be interpreted as proving that $E^{spin}_n$ is a simplicial $H$-space, as discussed in \S5.2.  Moreover, there is a  homotopy equivalence $E^{spin}_{n-1} \to \Omega E^{spin}_n$, as discussed in \S5.3.  Specifically, there is a functorial group isomorphism $[X, E^{spin}_{n-1}] \to [C^+X/ X, E^{spin}_n]$ defined by $s(w,p) = (sw, sp)$.  The key point here is the property $s (p \cup_{n-3} q) = sp \cup_{n-2} sq$ of the cochain suspension map given in Equation (4.1) in \S4.2.  Thus, the 2-stage Postnikov towers $E^{spin}_n$ form a simplicial 2-stage Postnikov spectrum $E^{spin}$.\\

\subsection {The Relation Between $E^{spin}$ and Spin Bordism}

We want to construct a natural isomorphism
$$G^{spin}_n(X)  \to {\rm Hom}(\widetilde{\Omega}^{spin}_n(X) / Image(\widetilde{\Omega}^{spin}_n(X^{n-2})),\ \R / \Z).$$
There is more than one way  to do this.  The entire discussion below applies to pairs $(X, Y)$ and relative bordism elements, but for simplicity we treat the case of a single space $X$. First, we can directly evaluate a pair $(w, p)$ on a reduced Spin bordism class $f\colon M^n \to X$.  A formula is
\begin{equation}
\langle (w, p), [M \xrightarrow{f} X] \rangle = \int_{[M]}\ f^*w\ + \  Q(f^*p) \in \R / \Z,
\end{equation}
where $Q$ is the canonical quadratic function on the Spin manifold $M^n$ studied in Chapter 2.  We gave two definitions of $Q$ in Chapter 2.  The arguments in Chapter 2 proving that the first (Kapustin) definition of $Q$ is well-defined and quadratic translate rather easily  to proving that this evaluation map is a well-defined group homomorphism. The image homomorphisms on bordism obviously vanish on $\widetilde \Omega_n^{spin}(X^{(n-2)})$, since there are no non-zero normalized cochains of degree $n-1$ and $n$ on the $n-2$ skeleton of $X$.\\

The second (Spin 1-manifold) definition of $Q$ is better suited for proving the evaluation map is an isomorphism.  This uses a filtration argument, more or less identical to the filtration argument given in ([1], \S6.3) in the case $n = 3$, which we will sketch three paragraphs below.\\

A completely different approach to an evaluation isomorphism, which will be discussed at the end of this section, is to use the fact that one knows in advance that the spectrum $E^{spin}$ that classifies the functors $G^{spin}_n(X)$ is the same as the spectrum that classifies the right hand bordism side of the desired isomorphism.  So, from stable homotopy theory, there is a natural group isomorphism, with certain properties that we will make explicit. But for now we continue studying the direct evaluation map (6.1).\\  

From the relations on representative pairs  $(w, p) \in G^{spin}_n(X)$,  we saw in \S3.1 that there is a short exact sequence
\begin{equation}
0 \to QH^n(X; \R /\Z) \to G^{spin}_n(X) \to SH^{n-1}(X; \Z/2) \to 0,
\end{equation}
where $$QH^n(X; \R /\Z) = H^n(X; \R / \Z)\ /\ Image(H^{n-2}(X; \Z/2) \xrightarrow{(1/2)Sq^2} H^n(X; \R / \Z))$$ and  $$SH^{n-1}(X; \Z/2) = Kernel( H^{n-1}(X; \Z/2) \xrightarrow{(1/2)Sq^2} H^{n+1}(X ; \R / \Z)).$$

On the Spin bordism side, there is a short exact sequence obtained as the Pontrjagin dual of the bordism sequence

\begin{equation}
0 \to Im(\tilde{\Omega}^{spin}_n(X^{n-1})) / Im(\tilde{\Omega}^{spin}_n(X^{n-2})) \to \tilde{\Omega}^{spin}_n(X)/Im(\tilde{\Omega}^{spin}_n(X^{n-2}))
\end{equation}
\begin{equation*}
 \to \tilde{\Omega}^{spin}_n(X) / Im(\tilde{\Omega}^{spin}_n(X^{n-1})) \to 0.
\end{equation*}
The sequence (6.2) will map to the Pontrjagin dual of the sequence (6.3), which will have arrow directions reversed, with the evaluation map in the center. From the Atiyah-Hirzebruch spectral sequence, the groups on the ends  in the dual of (6.3) are also $QH^n(X; \R/ Z)$ and $SH^{n-1}(X; \Z/2)$. The second definition of $Q$ can be used to show that the diagram of maps of short exact sequences commutes.  A Five Lemma argument completes the proof that the evaluation is an isomorphism.\\

Since both the theory represented by the spectrum $E^{spin}$ and the Pontrjagin dual of Spin bordism are cohomology theories, we want our evaluation comparison of the two theories to commute with the suspension isomorphisms in the two theories.  This is indeed the case and is explained by the cochain suspension map $s$ and transversality for bordism.   In the diagram below, the vertical arrows are the evaluation homomorphisms.  The top arrow is cochain suspension $s(w,p) = (sw, sp)$, with the homotopy theoretic interpretation of \S5.3.  Here we write $\Sigma X$ instead of $C^+ X/X $.
\begin{align}\label{diag1}
\begin{split}
[X, E_{n-1} ]\hspace{.3in}&\xrightarrow{s}\ \ \   [\Sigma X, E_n]   \\ 
\downarrow{e} \hspace{.5in}& \hspace{.5in}  \downarrow{e} \\
{\rm Hom}(\widetilde{\Omega}_{n-1}^{spin}(X),\ \R /\Z) & \xrightarrow{\sigma^*}  {\rm Hom}(\widetilde{\Omega}_n^{spin}(\Sigma X) ,\ \R /\Z) \\
\end{split}
\end{align}

The bottom arrow is the Pontrjagin dual of the geometric suspension isomorphism in bordism, $\sigma \colon \widetilde{\Omega}_n^{spin}(\Sigma X) \to \widetilde{\Omega}_{n-1}^{spin}(X)$.    Given $f \colon M^n \to \Sigma X$, transversality produces a collared submanifold  $(-1, 1) \times N \subset M$ and a map, $$(-1, 1) \times N^{n-1} \to (-1, 1) \times X \subset \Sigma X.$$  Then $\sigma(f) = f \vert_{0 \times N} \to 0 \times X$.  We need an orientation and Spin structure on $N$. The choice that makes the diagram commute is $N = (-1)^n \partial (M^+)$, where $M^+ = f^{-1}(C^+ X) \subset M$ and $\partial M^+$ is oriented with the outward (downward pointing) normal first.\footnote{If we had used the lower cone $C^- X$ to define cochain suspension in  \S4.1, we would need no sign if we take $N = \partial M^-$ oriented with the upward pointing normal first.} To prove the commutativity of the diagram, it is easiest to use the Spin 1-manifold version of the evaluation homomorphisms.  The sign $N = (-1)^n\partial (M^+)$ in the orientation of $N$ is needed to get Stokes Theorem to work out, $ \int_ {[M^+]}\ sw = \int_{[N]}\ w$, as part of the proof of commutativity of Diagram (6.4).\\

It is important to know that the natural evaluation isomorphism just defined is  essentially unique.  This follows from the following result.\\

CLAIM 6.1:  The functor $G^{spin}_n(X)$ admits no natural automorphisms other than the identity map and the inverse map of abelian groups.\\

To prove this, consider a natural automorphism $(w, p) \mapsto (w', p')$.  Then we must have $p' = p + dc$, since there are no natural automorphisms of $SH^{n-1}(Z; \Z/2)$ other than the identity.  So after multiplying by a relation $((1/2)Sq^2c, dc)$ we can assume $p = p'$.  Then $dw = dw' = (1/2)Sq^2p$.  Thus, $2w, 2w'$ and $w - w'$ are cocycles.  In the universal example, both $2w$ and $2w'$ represent generators of a (continuous) cohomology group isomorphic to $\Z$.  It follows that up to coboundaries $w$ and $w'$ are either equal or negatives.  This proves the Claim.\\

One then sees that the only other formula for a natural evaluation isomorphism would be to put a $-1$ in front of the $\int_{[M]}\ f^*w$ part of the evaluation formula (6.1), since $G^{spin}_n(X)$ admits no natural automorphisms other than $\pm Id$.\\

We will now sketch a second discussion of the comparison of $G^{spin}_n(X)$ with the Pontrjagin dual of Spin bordism of $X$. This second approach is the only one we understand in the Pin case.  The Postnikov tower of the spectrum $MSpin$ begins with a $K(\Z, 0)$, then a $K(Z/2, 1)$ with non-zero $k$-invariant  $Sq^2 \iota_0$.  By consideration of the Atiyah-Hirzebruch spectral sequence for the Pontrjagin dual cohomology theory, there is only one possibility for the classifying space of $$\rm{Hom}(\tilde{\Omega}^{spin}_n(X)/Im(\tilde{\Omega}^{spin}_n(X^{n-2})), \R / \Z),$$ and that is the 2-stage tower $E^{spin}_n$ we have been studying.  More precisely, the spectrum $E^{spin}$, as $n$ varies. Therefore, we know in advance that there is a diagram of  evaluation isomorphisms from sequence (6.2) to the dual of sequence (6.3), and this diagram of isomorphisms is functorial in $X$.  We also know in advance that we will have commutative suspension diagrams of homomorphisms (6.4), as $n$ varies, and these are also functorial in $X$.\\  

Now let's fix  $(X, Y) = (M^n, \partial M^n)$, a Spin manifold with boundary, and $f = Id \colon (M, \partial M) \to (M, \partial M)$.  As mentioned earlier, all the discussion above about evaluation maps applies to pairs $(X, Y)$.  By duality, there will be an associated evaluation map $G^{spin}_n(M, \partial M) \to \R / \Z$.    The group theory in \S3 implies the evaluation must have the form $$\langle (w, p), [(M, \partial M) \xrightarrow{Id} (M, \partial M)] \rangle = I(w) + \  Q(p) \in \R / \Z,$$ where $I$ is linear and $Q$ is quadratic.  From the map between diagram (6.2) and the dual of (6.3), we also know that $I(w) = \pm \int_{[M, \partial M]}\ w$.  Finally, given $u\colon (M^n, \partial M^n) \to (S^{n-1}, pt)$ with $p = u^*(z)$, the commutativity of the evaluation isomorphisms with $u^* \colon \Z/2 = G^{spin}_n(S^{n-1}, pt) \to G^{spin}(M^n, \partial M^n)$ implies that $Q(p)$ must be the Spin bordism class of $u$.\\

From the suspension diagrams (6.4), or more directly, the bordism class of $u$ is the same as the  Spin bordism class $[Z]$, where the  1-manifold $Z = u^{-1}(pt) \subset M$ is framed in $M$, hence has a Spin structure.  This is exactly the second definition of the canonical quadratic function $Q\colon Z^{n-1}(M, \partial M ; \Z/2) \to \Z/2 \subset \R / \Z$ on a Spin manifold given in \S2.3.\\

What we have proved at the end of this section is that the canonical quadratic function $Q$ on a Spin manifold constructed in \S2.2 and \S2.3 coincides with a canonical function  $Q$ that is necessarily quadratic, which arises from some stable homotopy theory and Pontrjagin duality. In particular, in \S2.3 we only knew that the second definition of $Q$ was quadratic because we could identify it with the first (Kapustin) definition in \S2.2.  But the proof here goes the other way.  We construct a function that must be quadratic, then we identify it with the second construction in \S2.3.  We will now turn to Pin bordism, and prove  by exploiting similar stable homotopy theory that the function $Q$ for Pin manifolds constructed in \S2.5 is necessarily quadratic.\\ 

\subsection{A Simplicial 2-Stage Postnikov Spectrum $E^{pin}$}

Consider the two stage Postnikov tower
$$ E^{pin}_n = K(\Z/2, n-1) \ltimes_{Sq^2p} K(\Z/2 , n).$$ 
Here, $p$ is the fundamental $\Z/2$ cocycle of degree $n-1$ and $Sq^2p = p \cup_{n-3} p$ is the standard cocycle representative of the cohomology operation $Sq^2$ in this dimension.   The space $E^{pin}_n$ represents the Pontrjagin dual  of a quotient of reduced $n$-dimensional Pin bordism.  Specifically, there is a natural isomorphism $$ G^{pin}_n(X) = [X, E^{pin}_n]  \simeq {\rm Hom}(\widetilde{\Omega}_n^{pin}(X) / Image\ \widetilde{\Omega}_n^{pin}(X^{(n-2)}),\ \R /\Z) $$ where $X^{(n-2)}$ denotes the $n-2$ skeleton of $X$.  We will explain this isomorphism in the next section.\\

Elements of $[X, E^{pin}_n]$ are represented by pairs $(w, p)$ with $dp = 0$ and $dw = p \cup_{n-3} p$.  The product will be given by
$$(w,p)(v,q) = (w+v + p \cup_{n-2} q,\  p+q).$$
The null-homotopic pairs turn out to be $(df +  Sq^2c,\ dc)$.   We now can follow exactly the discussion in \S6.1 about the spaces $E^{spin}_n$. It is only necessary to remove some $(1/2)$'s in front of cochains, since here we have only $\Z/2$ cochains, not $\R / \Z$ cochains.  Thus pairs modulo all null-homotopic pairs form an abelian group, and this group is exactly the group $G^{pin}_n(X) = [X, E^{pin}_n]$ studied in \S3.2.\\

The arguments above can be interpreted as proving that $E^{pin}_n$ is a simplicial $H$-space, as discussed in \S5.2.  Moreover, there is a  homotopy equivalence $E^{pin}_{n-1} \to \Omega E^{pin}_n$, as discussed in \S5.3.  Specifically, there is a functorial group isomorphism $[X, E^{pin}_{n-1}] \to [C^+X/ X, E^{pin}_n]$ defined by $s(w,p) = (sw, sp)$.  The key point here is the property $s (p \cup_{n-3} q) = sp \cup_{n-2} sq$ of the cochain suspension map given in Equation (4.1) in \S4.2.  Thus, the 2-stage Postnikov towers $E^{pin}_n$ form a simplicial 2-stage Postnikov spectrum $E^{pin}$.\\

\subsection{The Relation Between $E^{pin}$ and Pin Bordism}

The Postnikov tower spectrum constructed in the previous section classifies the cohomology theory Pontrjagin dual to $\tilde{\Omega}^{pin}_n(X)/Im(\tilde{\Omega}^{pin}_n(X^{n-2}))$, as $n$ varies.  This follows by considering the Atiyah-Hirzebruch spectral sequence for Pin bordism and its Pontrjagin dual,  because the $0^{th}$ and $1^{st}$ Pin bordism groups are both $\Z/2$, and the first $k$-invariant of $MPin$ is non-trivial.\\

The discussion in this section applies to pairs $(X, Y)$ of spaces, but to keep notation simpler we focus on the absolute case. From the relations on representative pairs  $(w, p) \in G^{pin}_n(X)$,  we saw in \S3.2 that there is a short exact sequence
\begin{equation}
0 \to QH^n(X; \Z / 2) \to G^{pin}_n(X) \to SH^{n-1}(X; \Z/2) \to 0,
\end{equation}
where $$QH^n(X; \Z /2) = H^n(X; \Z /2)\ /\ Image(H^{n-2}(X; \Z/2) \xrightarrow{Sq^2} H^n(X; \Z /2))$$ and  $$SH^{n-1}(X; \Z/2) = Kernel( H^{n-1}(X; \Z/2) \xrightarrow{Sq^2} H^{n+1}(X ; \Z /2)).$$

On the Pin bordism side, there is a short exact sequence obtained as the Pontrjagin dual of the bordism sequence

\begin{equation}
0 \to Im(\tilde{\Omega}^{pin}_n(X^{n-1})) / Im(\tilde{\Omega}^{pin}_n(X^{n-2})) \to \tilde{\Omega}^{pin}_n(X)/Im(\tilde{\Omega}^{pin}_n(X^{n-2}))
\end{equation}
\begin{equation*}
 \to \tilde{\Omega}^{pin}_n(X) / Im(\tilde{\Omega}^{pin}_n(X^{n-1})) \to 0.
\end{equation*}
Because the classifying spaces are the same, the sequence (6.5) will map naturally to the Pontrjagin dual of the sequence (6.6), which will have arrow directions reversed, with a group isomorphism in the center. From the Atiyah-Hirzebruch spectral sequence, the groups on the ends  in the dual of (6.6) are also $QH^n(X; \Z/2)$ and $SH^{n-1}(X; \Z/2)$.\\

Since both the theory represented by the spectrum $E^{pin}$ and the Pontrjagin dual of Pin bordism are cohomology theories, a natural evaluation isomorphism of the two theories will commute with the suspension isomorphisms in the two theories.  That is, we will have a commutative  diagram (6.7), where the vertical arrows are  evaluation homomorphisms.  The top arrow is cochain suspension $s(w,p) = (sw, sp)$, with the homotopy theoretic interpretation of \S3.3.  Here we write $\Sigma X$ instead of $C^+ X/X $.
\begin{align}\label{diag1}
\begin{split}
[X, E_{n} ]\hspace{.3in}&\xrightarrow{s}\ \ \   [\Sigma X, E_{n+1}]   \\ 
\downarrow{e} \hspace{.5in}& \hspace{.5in}  \downarrow{e} \\
{\rm Hom}(\widetilde{\Omega}_{n}^{pin}(X),\ \R /\Z) & \xrightarrow{\sigma^*}  {\rm Hom}(\widetilde{\Omega}_{n+1}^{pin}(\Sigma X) ,\ \R /\Z) \\
\end{split}
\end{align}

The bottom arrow is the Pontrjagin dual of the geometric suspension isomorphism in bordism, $\sigma \colon \widetilde{\Omega}_{n+1}^{pin}(\Sigma X) \to \widetilde{\Omega}_{n}^{pin}(X)$.    Given $f \colon M^{n+1} \to \Sigma X$, transversality produces a collared submanifold  $(-1, 1) \times N^n \subset M$ and a map, $$(-1, 1) \times N^{n} \to (-1, 1) \times X \subset \Sigma X.$$  Then $\sigma(f) = f \vert_{0 \times N} \to 0 \times X$.  The manifold $N$ has a Pin structure because it is framed in $M$.\\

It is important to know that the natural evaluation isomorphism
$$ G^{pin}_n(X) = [X, E^{pin}_n]  \simeq {\rm Hom}(\widetilde{\Omega}_n^{pin}(X) / Image\ \widetilde{\Omega}_n^{pin}(X^{(n-2)})\ \R /\Z) $$ is essentially unique.  This follows from the following result.\\

CLAIM 6.2:  The functor $G^{pin}_n(X)$ admits no natural automorphisms other than the identity map and the inverse map of abelian groups.\\ 

The proof is somewhat similar to the proof of Claim 6.1.  Express an automorphism as $(w, p) \mapsto (w', p')$.  Multiplying by a relation $(Sq^2c, dc)$, we can assume $p' = p$.  Now $dw = dw' = Sq^2 p$.  The only non-zero element in $H^n(E^{pin}_n ; \Z /2)$ is the class of $Sq^1p$.  Therefore $w' - w$ is cohomologous to either 0 or to $Sq^1 p$.  In the first case $(w', p) = (w, p)$ in  $G^{pin}_n(X)$.  In the second case $(w', p) = (w + Sq^1p, p) = (w, p)^{-1}$ in $G^{pin}_n(X)$.\\

\subsection{The Pin Bordism of $\R P^2$ and Other Surfaces}

We will now insert a computation of the reduced Pin bordism group $\tilde{\Omega}^{pin}_2(\R P^2)$, or more precisely its Pontrjagin dual.  For $X = \R P^2$  the short exact sequence (6.5) simplifies to
$$ 0 \to H^2(\R P^2; \Z/2) = \Z / 2 \to G^{pin}_2(\R P^2) \to \Z/2 = H^1(\R P^2; \Z/2) \to 0.$$
Pairs $(w, x)$ representing elements of $ G^{pin}_2(\R P^2)$, are given by arbitrary cocycles, since $dw = Sq^2x = 0$.  The product is given by $(w, x)(v, y) = (w+v+xy,\ x+y)$.  In particular $(0, x)^2 = (x^2, 0)$, so the extension is non-trivial and $ G^{pin}_2(\R P^2) \simeq \Z/4.$\\

If $(E^2, \partial E^2)$ is a connected  surface, we have the exact sequence $$ 0 \to H^2(E^2, \partial E^2; \Z/2) = \Z / 2 \to G^{pin}_2(E^2, \partial E^2) \to H^1(E^2, \partial E^2; \Z/2) \to 0.$$ The product formula in $G^{pin}_2(E^2, \partial E^2)$ is still $(w, x)(v, y) = (w+v+ xy, x+y)$.  Then $G^{pin}_2(E^2, \partial E^2)$ is a direct sum of $\Z/2$'s and $\Z/4$'s, with the number of $\Z/4$'s coinciding with the rank of the map $x^2 \colon H^1(E^2, \partial E^2; \Z/2) \to  \Z/2$.\\

We need a couple of other old facts.  The Pin bordism group $\Omega^{pin}_2(pt)$ is isomorphic to $\Z/8$.  This was established in the 1960's by Adams spectral sequence arguments.  $\R P^2$, with either Pin structure, is a generator.  In this Pin bordism group, the `negative' of $\R P^2$ with one Pin structure is $\R P^2$ with the other Pin structure\footnote{In general for a closed Pin surface $L^2$ the negative in the Pin bordism group is $\L^2$ with its Pin structure twisted by $w_1(L) \in H^1(L; \Z/2)$.}.  In particular, the identity maps $Id \colon \R P^2 \to \R P^2$, with the two different Pin structures on the domain, are {\em not} Pin bordant.\\

Let $K^2 $ denote the Mobius band, which we regard as an obvious subset of $\R P^2$. The computation  $G^{pin}(K^2, \partial K^2) \simeq \Z /4$ is the same as the computation for $\R P^2$.
Obviously there are two isomorphisms to $\Z / 4$, but fixing one is a somewhat delicate matter, as is fixing the evaluation isomorphism   $$e\colon G^{pin}_2(\R P^2) \simeq {\rm Hom}(\widetilde{\Omega}_2^{pin}(\R P^2),\ \R /\Z),$$ or the equivalent evaluation for $(K^2, \partial K^2)$.  Which Pin structure should be considered preferred on $\R P^2$ and  $K^2$?  Classically, a Pin structure on $M^n$ can be interpreted as a Spin structure on the bundle $\tau_M + det(\tau_M)$, where $\tau_M$ is the tangent bundle.  For a surface this is equivalent to a homotopy class of trivializations of  $\tau_M + det(\tau_M)$. As a choice of preferred Pin structure on $\R P^2$ and $K^2$, we will view  $\tau_{\R P^2} + det(\tau_{\R P^2})$ as an open submanifold of $\R P^3$, and give $\R P^3$ a Spin structure as the boundary of the tangent disk bundle of $S^2$ with its complex orientation. This is the same as the Spin structure given by the Lie group framing of $\R P^3 = SO(3)$.\\

One advantage of $K^2$ over $\R P^2$ is that $\R P^2$ is not a boundary, so the identity map is not a reduced bordism element.  We avoid fixing this by working with $K^2$ and relative bordism.  As the preferred generator  $(+1) \in \widetilde{\Omega}_2^{pin}(K^2, \partial K^2) \simeq \Z / 4$ we will  take $Id \colon (K^2, \partial K^2) \to (K^2, \partial K^2),$  where the domain $K^2$ is given the Pin structure described above.  The identity map with the other Pin strucuture on  the domain $K^2$ necessarily gives the other generator of $\Z/4$.  This follows from the comment above that the two Pin structures on $\R P^2$ are not Pin cobordant. \\

Now here are some useful facts that we will need in \S6.7 when we consider the effect of change of Pin structure on canonical quadratic functions.  Suppose $A^2 = I \times S^1$ is the annulus and suppose $u \colon (A^2, \partial A^2) \to (K^2, \partial K^2)$ is a  map of interval bundles over $S^1$.  Then a fiber $I$ maps to a fiber, so the induced maps on both relative $H_1$ with $\Z/2$ coefficients and relative $H^1$ with $\Z/2$ coefficients are isomorphisms.  Put another way, the induced map $u^* \colon \Z/4 = G^{pin}_2(K^2, \partial K^2) \to G^{pin}(A^2, \partial A^2) = \Z/2 \oplus \Z/2$ is non-zero.  In fact, the map is projection onto the $\Z /2$ direct summand spanned by $(0, x)$.\\

There are two Pin structures on the annulus, which correspond by evaluation to the two (linear) maps $G^{pin}(A^2, \partial A^2) \to \Z/2$ that factor through $H^1(A^2, \partial A^2 ; \Z/2)$.  Therefore, given $u \colon (A^2, \partial A^2) \to (K^2, \partial K^2)$ as in the above paragraph, then with the non-trivial Pin structure on $A^2$ any such  map $u$ represents $2 \in \Z/4$ in the Pin bordism group of $(K^2, \partial K^2)$.  With the trivial Pin structure on $A^2$, the map represents 0 in Pin bordism.\\

Similarly, suppose $u \colon (K^2, \partial K^2) \to (K^2, \partial K^2)$ is a map of interval bundles over $S^1$, which necessarily has odd degree on the core circles.  Then the induced map $u^*$ from $G^{pin}_2(K^2, \partial K^2) = \Z/4$ to itself is the identity map. This holds because the generator $(0, x) \in G^{pin}_2(K^2, \partial K^2)$ pulls back to itself.  We conclude that with the preferred Pin structure on the domain $K^2$, any such map $u$ represents $+1 \in \Z/4$ in the Pin bordism group of $(K^2, \partial K^2)$, and with the opposite Pin structure on the domain $K^2$, such a map represents $-1 \in \Z/4$.

\subsection{The Canonical Quadratic Function on Pin Manifolds} 

Now let's fix  $(M^n, \partial M^n)$ with a Pin structure, and consider $Id \colon (M, \partial M) \to (M, \partial M)$ as a relative Pin bordism element.  By duality, there will be an associated evaluation map $G^{pin}_n(M, \partial M) \to \R / \Z$.    The group theory in \S3 implies that an evaluation must have the form $$\langle (w, p), [(M, \partial M) \xrightarrow{Id} (M, \partial M)] \rangle = I(w) + \  Q(p) \in \R / \Z,$$ where $I$ is linear and $Q$ is quadratic.  From the map between diagram (6.5) and the dual of diagram (6.6), we  know that $I(w) =  (1/2)\int_{[M, \partial M]}\ w$.  In the Spin case, the sign ambiguity of the corresponding evaluation map implicit in Claim 6.1 was pinned down by the choice that $(w, 0)$ should evaluate as $\int_{[M, \partial M]}\ w \in \R / \Z$.  In the Pin case, the integral has order 2, so the sign ambiguity implicit in Claim 6.2 must be resolved by explaining how to choose  $Q(p) \in \Z / 4 \subset \R / \Z$.\\

The universal sign ambiguity is pinned down by the decisions made for $M = \R P^2$ or $K^2$.  In those cases there are obviously two quadratic functions $q(x) = \pm 1 \in \Z /4$, where $(0,x) $ generates $G^{pin}_2$.  With the preferred Pin structure on $\R P^2$ and $K^2$ defined in the previous section, we will choose as the preferred quadratic function $q(x) = +1$.  This is the same as choosing the evaluation map for $\R P^2$ to satisfy $$\langle (0, x), [ \R P^2 \xrightarrow{Id} \R P^2] \rangle = +1/4 \in \R / \Z,$$ where the domain $\R P^2$ is given the preferred Pin structure.\\

We first discuss  the canonical quadratic function on surfaces, beginning with closed surfaces.  Recall in dimension two that quadratic functions are defined on cohomology groups.  Given any $p \in H^1(E^2; \Z/2)$, there will be maps $u\colon E^2 \to \R P^2$ with $u^*(x) = p$.  Then the map $u^* \colon G^{pin}_2(\R P^2) \to G^{pin}_2(E^2)$ with $u^*(0, x) = (0, p)$ will commute with the natural evaluation maps to the duals of Pin bordism.  That is, there is a commutative diagram
\begin{align}
\begin{split}
G^{pin}_2(E^2)\hspace{.3in}&\xleftarrow{u^*}\ \ \   G^{pin}_2(\R P^2)  \\ 
\downarrow{ev} \hspace{.5in}& \hspace{.8in}  \downarrow{ev} \\
{\rm Hom}(\widetilde{\Omega}_{2}^{pin}(E^2),\ \R /\Z) & \xleftarrow{Du_*}  {\rm Hom}(\widetilde{\Omega}_{2}^{pin}(\R P^2) ,\ \R /\Z) \\
\downarrow{ev(E^2)} \hspace{.2in} &  \hspace{.8in} \downarrow{ev(\R P^2)}\\
\R / \Z \hspace{.6in} &   =  \hspace{.5in}   \R / \Z \\
\end{split}
\end{align}

We then see the following key fact.  Given our normalization conventions, the Pin bordism class $[u \colon E^2 \to \R P^2] \in \Z/4 \subset \R /\Z$ coincides with the value of the canonical quadratic function $Q(p) = Q(u^*x) \in \Z / 4 \subset \R / \Z.$\\

The case of Pin surfaces with boundary is essentially identical to the above, exploiting maps $u \colon (E^2, \partial E^2) \to (K^2, \partial K^2) $, with $u^*(x) = p$.  The Pin bordism class of $u$ coincides with the value of the quadratic function $Q(p) = Q(u^*x)$.  Both take well-defined values in $\Z/4$, given the normalization conventions for $K^2$.\\

In higher dimensions, suppose given $u\colon (M^n, \partial M^n) \to (\Sigma^{n-2} \R P^2, pt)$, with $p = u^*(\sigma^{n-2}x)$ where $x \in Z^1(\R P^2; \Z/2)$ is non-zero in cohomology.  The commutativity of the natural evaluation maps with $u^* \colon G^{pin}_n (\Sigma^{n-2} \R P^2, pt) \to G^{pin}_n(M^n, \partial M^n)$, and the commutativity of suspension diagrams (6.7), imply that $Q(p) \in \Z / 4 \subset \R / \Z$ must be the Pin bordism class of $u$, which is the same as the  reduced Pin bordism class  $[u \colon E^2 \to \R P^2] \in \Z/4$ where the Pin surface $E^2 \subset M$ is the framed inverse image of $\R P^2$ in $M$ under the map $u$.\\

This proves that our formula for $Q(p)$ back in \S2.5, namely $$Q(p) = (1/4)[u \colon E^2 \to  \R P^2] \in \R / \Z$$ in the case that  $p = u^*(\sigma^{n-2}x)$), does indeed  define\footnote{Recall that the general cocycle $p'= p + dc$, where $p$ is of this special form.  Then the quadratic  properties of $Q$ give the formula for $Q(p+dc)$.} a quadratic function $Z^{n-1}(M, \partial M; \Z/2) \to \R / \Z$.  This is the quadratic function that we will declare to be the canonical quadratic function on a Pin manifold $M^n$.  Of course it is only really canonical as a $\Z/4 \subset  \R / \Z$ valued function after we make the normalization choices described above for $\R P^2$ and $K^2$.\\

We make a few last comments about the definition of the canonical quadratic function.  Historically, and in great detail in reference [4], the canonical $\Z / 4$-valued quadratic function on a surface was defined by interpreting a Pin structure as a trivialization of the bundle $\tau_{E} + det(\tau_{E})$, then counting `right hand twists' of a band around a circle in $E$.  The count is in $(1/2)\Z / 2 \Z \simeq \Z/4$.  We have not needed to identify our canonical quadratic function with that one. But they are the same.  Resolution of the sign ambiguity occurs as the decision to count right hand twists, rather than left hand twists.  Also, one must decide at some point which is the `true' Pin $\R P^2$, to settle the ambiguity $q(x) = \pm 1$. Proof that in the surface case the twist counting yields quadratic functions defined on cohomology and refining the cup product pairing  is done geometrically, since the cup product can be viewed as self-intersection, which can be understood.  Not so much with the higher cup products in our $n$-dimensional  quadratic functions defined on cocycles.\\

On a simplicial $n$-manifold $M^n$, one can draw curves $Z$ dual to an $n-1$-cocycle $p$ by putting points on $n-1$ simplices on which $p$ evaluates non-trivially and connecting with arcs in $n$ simplices.  Identify the normal bundles around component circles of $Z$, and use a Pin structure to make a twist count of some kind in framed 3 dimensional tubes around these circles.  The catch is, what does quadratic mean, as a function of $p$, and how do you prove it?  Curves for different cocycles $p, q$ get close together when $p$ and $q$ both evaluate non-trivially on the same $n-1$ simplex, or even on faces of the same $n$-simplex, and this is exactly when $p+q$ gets interesting. Our notion of quadratic function is defined in terms of the black box $\cup_{n-2}$ operation.  The cochain formulas are very complicated. We see no way to handle this in a direct geometric manner.  Therefore, we use the method of Chapter 2, or the even more abstract method of Chapter 6, to get quadratic functions, which we know must be identical to some twist count construction.\\

\subsection{Change of Pin Structure on a Manifold}

In order to complete the proof of our original Claim 1.1 in the Pin case, we need to understand how the canonical quadratic function $Q$ on $M$ is related to the canonical quadratic function on $M_a$ if the classical Pin structure on $M$ is changed to $M_a$, by the action of $a \in H^1(M; \Z/2)$ on Pin structures.  The result we want is that the canonical quadratic function on $M_a$ is indeed $Q_a(p) = Q(p) + 2\langle ap, [M, \partial M] \rangle$.   Note that with $p = u^*(\sigma^{n-2}x)$ we have $\langle ap, [M, \partial M] \rangle = \langle i^*(a)u^*(x), [E^2]\rangle$, where $i \colon E^2 \subset M$ is framed  by the map $u\colon M \to \Sigma^{n-2}\R P^2$. Thus we are reduced to the following question about surfaces and a reduced Pin bordism element $u\colon E^2 \to \R P^2$, with $u^*(x) = p \in H^1(E^2; \Z / 2)$.  How does the Pin bordism class of $u$ change if the Pin structure on $\Sigma$ is changed by $a \in H^1(E^2; \Z/2)$? The result we want is that it changes by $2\langle ap, [E^2] \rangle \in \Z /4$.\\

First, we have the Mobius strip $K^2 \subset \R P^2$ and we make $u$ transverse to the non-trivial circle in $\R P^2$, so that its inverse image is a union of circles, with a tubular neighborhood $E_0$ mapping to $K^2$.  We then replace $u$ with $u_0\colon\ (E_0,\ \partial E_0) \to (K^2,\ \partial K^2)$.  Then $E_0$ is a union of annuli and Mobius strips, with Pin structures from $E$, whose core circles map to the core circle in $K^2$ with even or odd degrees, respectively. Which elements in the Pin bordism of $(K^2, \partial K^2)$, which we have identified with $\Z/4$,  do these annuli and Mobius strips represent?  The answer was given in the final three paragraphs of \S6.5.   A Mobius strip mapping with odd degree  to $K^2$ represents  $\pm 1 \in \Z /4$, where the sign compares the Pin structure on the Mobius  strip in $E_0$ with the chosen preferred Pin structure on $K^2$.  An annulus mapping with even degree to $K^2$ represents $2 \in \Z/4$ if the Pin structure on the annulus in $E_0$ is non-trivial, and represents $0 \in \Z/4$ if the Pin structure from $E_0$ is trivial. \\

OK.  Now we are ready to see what happens to $u \colon (E_0, \partial E_0) \to (K^2, \partial K^2)$ if we change the Pin structure on $E_0$ by $a \in H^1(E_0; \Z /2)$.  We interpret the Pin structure change as twisting the normal framings of  circles $Z$ in $\tau_{E} + det(\tau_E)$ by the number $\langle a, [Z] \rangle \in \Z /2$.  This coincides with $\langle ap_z, [E_0, \partial E_0] \rangle \in \Z/2$ if $p_z \in H^1(E_0, \partial E_0; \Z /2)$ is dual to $[Z]$.  Reviewing the discussion in the paragraph above, this gives exactly what we want.  Given any $p = u^*(x) \in H^1(E_0, \partial E_0; \Z /2)$, where $u \colon (E_0, \partial E_0) \to (K^2, \partial K^2)$,  the dual of $p$ is a union of circles, with neighborhoods mapping to $K^2$.  The analysis above shows that the relative Pin bordism class represented by a neighborhood of one  circle $Z$  changes by $2\langle a, [Z]\rangle \in \Z /4$.  So the global change on $(E_0, \partial E_0) \to (K^2, \partial K^2)$ is exactly the sum of these, or $2\langle ap, [E_0, \partial E_0] \rangle.$\\

\section{Some Manipulations with Quadratic Functions }

In this Chapter we make several direct constructions with quadratic functions.  Some of these are analogues of obvious constructions with classical Pin structures.  For example, a Pin structure on a manifold induces a Pin structure on its boundary and on codimension zero submanifolds.  Also, there is an obvious notion of cobordism between Pin manifolds. But here we want to make these constructions directly with quadratic functions, without reference to  the main Claim 1.1(P), or its proof in the Chapters above.  Other constructions in this section are perhaps not so familiar with classical Pin structures.\\

\subsection{Quadratic Functions on Boundaries}

Here is how we construct quadratic functions on boundaries. We are assuming that on each simplex of $M$ the vertices on $\partial M$ precede the vertices not on $\partial M$.   We then have a collapse map in the ordered simplicial category $t \colon (M, \partial M) \to (C^+\partial M, \partial M)$, extending the identity on $\partial M$ and mapping all other vertices of $M$ to the cone point. Note $t$ is an isomorphism on each $n$-simplex of $M$ that meets the boundary in an $(n-1)$-simplex, and maps all other $n$-simplices  degenerately.\\

We then have cochain maps $$t^* s\colon Z^{n-2}(\partial M ;\Z/2) \to Z^{n-1}(C^+\partial M, \partial M; \Z/2) \to Z^{n-1}(M, \partial M ;\Z/2),$$ where $s$ is the cochain suspension studied in \S4.1.  For a cocycle $w \in Z^j(\partial M; \Z/2)$ of any degree, it is  useful to observe $t^*s(w) = d\hat{w}$, where $\hat{w} \in C^j(M; \Z/2)$ extends $w$ by 0 on all simplices not in $\partial M$.  So $t^*sw$ is an explicit relative cocycle representative for $\delta\bar{w} \in H^{j+1}(M, \partial M; \Z/2).$  Thus  for any  $u \in C^{n-1}(\partial M; \Z/2) = Z^{n-1}(\partial M; \Z/2)$ we  have $$\int_{[M, \partial  M]}\ t^*su\ =\ \int_{[C^+M, \partial M]}\ su\ =\    \int_{\partial M}\ u.$$
We now compose with a quadratic function $Q$ on $(M, \partial M)$, $$\partial Q = Q\circ t^*s \colon Z^{n-2}(\partial M; \Z/2) \to Z^{n-1}(M, \partial M; \Z/2) \to \Z/2.$$
We claim $\partial Q$ is a quadratic function on $\partial M$.  First, $$\partial Q(p+q) = \partial Q(p) + \partial Q(q) + \int_{[C^+ \partial M, \partial M]} sp \cup_{n-2} sq.$$ But $sp \cup_{n-2} sq = s(p \cup_{n-3}q)$,  by property (4.1) of the cochain suspension $s$ in \S4.2.  The  integral term then equals $\int_{\partial M} p \cup_{n-3} q$.  Secondly, if $c \in C^{n-3}(\partial M; \Z/2)$ then $$\partial Q(dc) = Q(dt^*(sc))  =\int_{[M, \partial M]}\ Sq^2(t^*sc) =  \int_{[C^+ \partial M, \partial M]}\ Sq^2 sc.$$ Again, by properties of $s$, the  integral equals $  \int_{\partial M} Sq^2 c.$  This completes the proof that $\partial Q$ is a quadratic function on $\partial M$.\\ 

Our proof of  Claim 1.1(P) constructed a canonical quadratic function on a Pin manifold, using classical facts about Pin structures.  So there is a canonical quadratic function on both a Pin manifold $M$ and on its boundary $\partial M$. If one examines the proof of Claim 1.1(P), it can be seen that the canonical quadratic functions constructed on $M$ and $\partial M$ are indeed related by the cochain suspension boundary construction on quadratic functions given here.\\

Although our construction of boundary quadratic functions uses only  direct operations with cochains, one can look at it from a more homotopy theoretic viewpoint.  The map  $M / \partial M \to C^+ \partial M / \partial M$ and the construction in \S5.3 yields a composition
$$[\partial M, E_{n-1}] \to  [C^+ \partial M / \partial M, E_n] \to [M / \partial M, E_n].$$ A quadratic function on $M$ is named by a  homomorphism $[M / \partial M, E_n] \to \R /\Z$.  Apply the functor $\rm{Hom}( * , \R / \Z)$ to the above composition and follow a quadratic function on $M$.  At the other end, one gets a homomorphism $[\partial M, E_{n-1}] \to \R / \Z$, which is equivalent to a quadratic function on  $\partial M$.\\

\subsection{Some Functorial Properties}

Next we point out some direct manipulations that amount to certain  functorial properties of quadratic functions.  The motivation for this is to indicate how one can develop some self contained theory of pairs $(M, Q)$ consisting of manifolds with a quadratic function. For example, if $f\colon (M', \partial M') \to (M, \partial M)$ is an order preserving simplicial map between manifolds of the same dimension such that $H^n(M, \partial M; \Z /2) \to H^n(M', \partial M' ; \Z/2)$ is surjective, and if $Q'$ is a quadratic function on $M'$ then it is obvious that  the composition $$Q = Q'  f^*\colon Z^{n-1}(M, \partial M; \Z /2) \to Z^{n-1}(M', \partial M'; \Z /2)  \to \Z / 2$$ is a quadratic function on $M$.  This is a push-forward construction.\\

A specific example is provided by the canonical order preserving map $M' \to M$, where $M'$ is the barycentric subdivision of $M$, with its natural ordered structure.  But also in this case, since this map induces a homology isomorphism,  a quadratic function $Q$ on $M$ induces a quadratic function $Q'$ on $M'$. The reason is, any $Q' $ is determined by the values of $Q'$ on cocycles $\{p_j'\}$ representing a homology basis of $H^{n-1}(M', \partial M' ; \Z/2)$, and such cocycles can be taken to be the image of cocycles $\{p_j\}$ on $(M, \partial M).$  Then set $Q'(p_j') = Q(p_j)$. This is a pull-back construction. The two constructions $Q \leftrightarrow Q'$ are clearly inverses of each other. In fact, these remarks show that any order preserving simplicial map between two $n$-manifolds inducing a cohomology isomorphism with $Z/2$ coefficients in degrees $n-2, n-1, $ and $n$ will induce a bijection between quadratic functions on the two manifolds.\\

\subsection {Low Dimensions}

It is amusing to look at low dimensions.  Logically speaking, when $n = 0$ the Claim 1.1(P) asserts every  0-manifold admits a unique Pin structure.  When $n = 1$, the Claim asserts that Pin structures correspond bijectively with linear functions $H^0(M, \partial M; \Z/2) = H^0(M';\Z/2) \to \Z/2$, where $M'$ is the union of the closed components of $M$.  (Such linear functions indeed correspond canonically with elements of $H^1(M; \Z/2)$.) In particular, when $M = S^1$ there is a canonical trivial `quadratic function' and a canonical non-trivial `quadratic function', meaning the zero and the non-zero homomorphisms $H^0(S^1, \Z/2) \to \Z/2$. \\

When $n = 2$ and $\partial M = \emptyset$,  Claim 1.1(P) is well-known, [4], and asserts that there is a canonical correspondence between Pin structures and quadratic functions $q\colon H^1(M^2; \Z/2) \to \Z/4$ refining the cup product pairing. If $\partial M \not= \emptyset$, the analysis gets a bit tricky. Pin structures still correspond canonically to quadratic refinements  of the cup product pairing for the pair $(M, \partial M)$, but the cup product pairing can be degenerate.  The restriction of quadratic functions induced by inclusions $\partial M^2 \subset M^2$ when $M^2$ is connected  reveals that the number of non-trivial homomorphisms induced on boundary circles must be even.  This will be explained in the next paragraph. Thus only the trivial  structure on $S^1$ bounds. \\

We will  look at specific examples of restricting quadratic functions to boundaries and other submanifolds.  We first examine the above boundary construction when $M^2$ is a connected 2-manifold.  Consider a cocycle $x_0 \in Z^0(\partial M; \Z/2)$  that is identically 1 on vertices of one component $Z_0$ of $\partial M$ and identically 0 on other components.  If $Q$ is a quadratic function on $(M, \partial M)$, then $ \partial Q(x_0) = 1\ \text{or}\ 0$  records whether the induced  structure  on $Z_0$ is the non-trivial or the trivial  structure (in our sense).  Let $x$ be the sum of these $x_0$ over the non-trivial components, and let $y$ be the sum over the trivial components. Then $t^*sx + t^*sy = dc$, where $c \in C^0(M, \partial M; \Z/2)$ is the relative cocycle that is 1 on all vertices of $M$ not in $\partial M$.  Then, because of the low dimension, $0 = Q(dc) = \partial Q(x) + \partial Q(y)$.  By definition, $\partial Q(y)$ is a sum of 0's.  Therefore, $\partial Q(x) = 0$, which says that the number of non-trivial boundary components must be even.  It is not hard to construct a quadratic function on $[0,1] \times S^1$ with boundary two copies of $S^1$ with non-trivial  structure.\\  

We next look at the boundary construction when $M^3$ is a connected 3-manifold. Then $\partial Q$ is a quadratic function on $H^1(\partial M; \Z/2)$.  Half the elements of this group are represented by cocycles $x = i^*z$ with  $z \in Z^1(M; \Z/2)$.  Then  $c = \hat{x} +z \in C^1(M, \partial M; \Z/2)$, where $\hat x$ extends $x$ by 0, and $dc = d\hat{x} = t^*sx$.  Then $Sq^2c = cdc = cd\hat x$, hence $\partial Q(x) = Q(dc) = \int_{[M, \partial M]} cd\hat{x}$.  But one sees that $cd\hat{x}$ is 0 on all 3-simplices since $c$ is 0 on all 1-simplices in $\partial M$ and $d\hat{x}$ is 0 on all 2-simplices disjoint from $\partial M$. Thus $\partial Q$ vanishes on a Lagrangian subspace of $H^1(\partial M; \Z/2)$, which, of course, implies its Arf invariant is 0.\\ 

\subsection {Codimension Zero Submanifolds}

In addition to restricting quadratic function to boundaries, we can also restrict to codimension 0 submanifolds $V$.  First assume $M$ is closed and $V \subset M$.  Cocycles $x \in Z^{n-1}(V, \partial V; \Z/2)$ can be extended by 0 to cocycles $\tilde{x} \in Z^{n-1}(M; \Z/2)$. We then define $Q_V(x) = Q(\tilde{x})$.  The desired two conditions are easy to check.\\

If $\partial M \not= \emptyset$ and if $V \subset M$ is neatly embedded, with $\partial V = \partial_0 V \cup \partial_1 V$, where $\partial_0 V = \partial V \cap \partial M$, then again we can extend relative cocycles $x$ on $(V, \partial V)$ by 0 to relative cocycles $\tilde{x}$ on $(M, \partial M)$, and define    $Q_V(x) = Q(\tilde{x})$. Iterated boundary and codimension 0 constructions are consistent under the various inclusions, $\partial_0 V \subset \partial V \subset V \subset M$ and $\partial_0 V \subset \partial M \subset M$.\\

Another aspect of the boundary construction and the codimension 0 construction is that  if $V' \subset M$ is the complementary codimension 0 submanifold to $V \subset M$, then we have $\partial_1 V = \partial_1 V'$ and $\partial_0 V \cup \partial_0 V' = \partial M$. We can restrict a quadratic function on $M$ to both $V$ and $V'$. Then further restrict to $\partial V' \supset \partial_1 V'$ and $\partial V \supset \partial_1 V$.  In the oriented setting, the orientations on $\partial_1 V$ and $\partial_1 V'$ are opposite, but it turns out the two induced quadratic functions  are the same. This is a fairly tricky computation. \\ 

\subsection{Cobordism and Homotopy}

Since we can restrict quadratic functions to boundaries  it is easy to define the notion of a cobordism $(W, Q)$ between $(M_i, Q_i),\ i = 0, 1$, where the $M_i$ are closed manifolds with  quadratic functions $Q_i$.  Cobordism of quadratic functions makes sense in both the oriented and the unoriented cases. In fact, since we can also restrict quadratic functions to codimension 0 submanifolds, we can define  relative cobordisms $(W, Q)$ between $(M_i, Q_i)$ when the boundaries of the $M_i$ need not be empty.  The construction will be in terms of the usual notion of a relative cobordism $(W, \partial W)$ between manifolds with boundary, together with a quadratic function $Q$ on $W$ that restricts appropriately to a given quadratic function on  the codimension 0 submanifold $M_0 \sqcup M_1 \subset \partial W$. \\

One can also go back and forth between quadratic functions on $M$ and quadratic functions on the specific (relative, if $\partial M \not= \emptyset$) cobordism $I \times M$.  Starting with a quadratic function $\widehat{Q}$ on $ I \times M$, restrict to the boundary and obtain quadratic functions $Q_i$ on $M_i = \{i\} \times M$.  Conversely, starting with a quadratic function $Q_0$ on $M_0$ exploit the homology isomorphism with a dimension shift induced by cochain suspension $s_0$ at the 0 end. Set $\widehat{Q}_0(s_0(x)) = Q_0(x)$ on representatives of a cohomology basis. This extends to a  quadratic function $\widehat Q_0$ on $I \times M$ that restricts to the original quadratic function $Q_0$ on $M_0$.  One sees that the two constructions $\widehat Q \leftrightarrow Q_0$ are bijective correspondences and inverses, since $H^1(I \times M, \Z/2) \simeq H^1(M; \Z/2)$. Now restrict $\widehat{Q}_0$ to $M_1$.  Call this restriction $Q_1$.  Go the other way.  Extend $Q_1$ to $\widehat{Q}_1$, using cochain suspension $s_1$ at the 1 end.  Since the restrictions of $\widehat{Q}_0$ and $\widehat{Q}_1$ agree on $M_1$, they must be identical. You are back where you started, $\widehat{Q}_1 = \widehat{Q}_0$.  This seems quite difficult to see directly unless the ordered simplicial structure on $I \times M$ is  very simple, something like a product triangulation extending the same simplicial structure on the two ends. \\

It is now pretty easy to `free' the notion of quadratic function from a choice of ordered simplicial structure on $M$, and to formulate a kind of functorial homotopy invariance. There are a couple of  ways to do this.  One can exploit iterated barycentric subdivisions, and declare quadratic functions $Q_0$ and $ Q_1$ on two ordered simplicial structures $M_0$ and $M_1$ on the same manifold $M$ to be equivalent if there is some ordered structure on $M$ and order preserving maps $f_i \colon (M, \partial M) \to (M_i, \partial M_i)$ homotopic to the identity so that the pull-backs to $M$ of the two quadratic functions $Q_i$ are identical. Or, perhaps more directly, just take arbitrary quadratic functions on ordered triangulations of $I \times M$ and declare the restrictions to the ends to be equivalent.  Related statements are that the push-forward  correspondences for two odd degree topologically homotopic ordered simplicial maps $f_0, f_1 \colon (M', \partial M') \to (M, \partial M)$ are identical, and similarly for the two pull-back correspondences when defined.\\

One can iterate the constructions with boundaries of $n$-manifolds and codimension 0 submanifolds  to deal with  some embedded submanifolds  of lower dimensions. For example, by means of such iterations  one can see that if $\Delta^{n-1} \times Z \subset M^n - \partial M^n$ is a framed circle in an oriented manifold dual in a suitable sense to a cocycle $p \in Z^{n-1}(M, \partial M;\Z/2)$ then $Q(p) = [Z] \in \Z/2$, where $Q$ is a quadratic function  on $M$ and $[Z]$ records whether the structure on $Z$ (in our sense) induced by iteration is trivial or non-trivial.  In the proof of the Claim 1.1 in the Spin case, we encountered the classical notion of Spin structure on a 1-manifold and that proof certainly shows that the two notions of Spin structure on a circle coincide.  A somewhat more elaborate discussion like this would apply to a  framed surface in an arbitrary manifold with a quadratic function.\\

An interesting point is that one can define a  bordism-like homology theory in both the oriented and non-oriented cases, by looking at bordism classes of maps $(M, \partial M ; Q) \to (X, Y)$, where $M $ is a manifold with a quadratic function $Q$.  The constructions in this Chapter yield proofs of the various homology theory axioms, namely homotopy functor, long exact sequence, excision. After the fact, one knows that these homology theories, combined with $M$ smooth, must coincide with Spin and Pin bordism.   Nonetheless, because the constructions are directly combinatorial in terms of simplicial structure, and we can cut a manifold into pieces and look at quadratic functions on the pieces, it seems like there is some connection here with extended TQFT theories.

\section* {\bf References}

1. Greg Brumfiel and John Morgan, The Pontrjagin Dual of 3-Dimensional Spin Bordism, arXiv.org$>$math.AT$>$arXiv:1612.02860v2.\\

\noindent2.  Greg Brumfiel and John Morgan, The Pontrjagin dual of 4-Dimensional Spin Bordism, arXiv.org$>$math.GT$>$:1803.08147v1.\\

\noindent3. D. Gaiotto and A. Kapustin,  Spin TQFTs and fermonic phases of matter,
arXiv.org$>$cond-mat$>$arXiv:1505.05856v2.\\ 

\noindent4. R. C. Kirby and L. R. Taylor, Pin structures on low-dimensional manifolds. Geometry of Low-dimensional Manifolds, 2, London Math. Soc. Lecture Note Ser., 151,
Cambridge Univ. Press (1990) 177-242. 

\end{document}